\documentclass[12pt]{article}

\usepackage{t-angles}
\usepackage[all]{xy}
\usepackage{amsmath,amsthm, amscd, amssymb, amsfonts}

\usepackage{hhline}

\usepackage{amsmath}
\usepackage{t-angles}

\title{
Duality Theorem and Hom Functor in Braided Tensor Categories}
\author{
Yange Xu $^{a, b}$, \ \ Shouchuan Zhang$^a$, \ \ Jing Cheng $^a$
 \ \  \ \\ a: Department  of Mathematics, Hunan University,
 \\ Changsha  410082, \
 P.R. China.\\
 b: Department  of Mathematics,  Pingdingshan University, \\ Pingdingshan
467000, \
 P.R. China.
  }
\date{}

\begin{document}
% \tableofcontents
\newtheorem{Theorem}{\quad Theorem}[section]
\newtheorem{Proposition}[Theorem]{\quad Proposition}
\newtheorem{Definition}[Theorem]{\quad Definition}
\newtheorem{Corollary}[Theorem]{\quad Corollary}
\newtheorem{Lemma}[Theorem]{\quad Lemma}
\newtheorem{Example}[Theorem]{\quad Example}
\maketitle \addtocounter{section}{-1}

\begin {abstract} The
Blatter-Montgomery duality theorem  is generalized into braided
tensor categories. It is shown  that  $Hom(V,W)$ is a  braided
Yetter-Drinfeld module for any two   braided Yetter-Drinfeld modules
$V$ and $W$.

 \end {abstract}
\section {Introduction} The duality theorems play an important role in actions of Hopf algebras (see \cite {Mo93}).
In \cite {BM85} and \cite {Mo93},  Blattner and Montgomery proved
the following  duality theorem for an ordinary Hopf algebra $H$ and
some Hopf subalgebra $U$ of $H^{*}$ :
\begin {eqnarray*} \label {e1}
(R \# H)\# U   \cong R \otimes (H \# U)    \hbox { \ \ \ as
algebras, \  }
 \end {eqnarray*}
 where $R$ is a $U$-comodule algebra. The dual theorems for co-Frobenius Hopf algebra $H$,
$$   (R \# H)\# H^{* rat}   \cong M_H^f (R) \ \hbox { and }
(R \# H^{* rat })\# H   \cong M_H^f (R)\hbox { \ \ \  as  } k \hbox
{-algebras  } $$ were proved in \cite {DZ99} and \cite {DNR01} (see
\cite [Corollary 6.5.6 and Theorem 6.5.11  ] {DNR01}). Braided
tensor categories become more and more important. They have been
applied in conformal field, vertex operator algebras, isotopy
invariants of links (see\cite {Hu05,HK04, BK01} and \cite {He91,
Ka97, Ra94}), One of the authors in \cite {Zh03} generalized the
duality theorem to the braided case, i.e., for a finite Hopf algebra
$H$ with $C_{H,H} = C_{H,H}^{-1},$
  \begin {eqnarray*}
(R \# H)\# H^{\hat *}   \cong R \otimes (H \bar
  \otimes H^{\hat *} )  \hbox { \ \ \ as
algebras in } {\cal C}.
 \end {eqnarray*} The Blattner-Montgomery  duality theorem was also generalized into Hopf algebras over commutative rings
\cite {ATL01}. $Hom$ functor also has extensive use in homological
algebra and representation theory.

We know that  $H$ is an infinite braided Hopf algebra if it has no
left duals (See \cite {Ta99}). In this paper we generalize the above
results to infinite braided Hopf algebras. In section 1,  we
introduce quasi-dual $H^d$ of $H$ and prove the duality theorem in a
braided tensor category ${\cal D}$; In section 2, we prove that if
$V, W$ are in $^{B}_{B}\mathcal{YD}(\mathcal{C})$, then $Hom(V,W)$
is also in $^{B}_{B}\mathcal{YD}(\mathcal{C})$; In section 3, we
concentrate on the Yetter-Drinfeld module category $^B_B {\cal YD}$.

 {\bf Some notations}. Let ($\mathcal{D},\otimes,I,C $) be a braided tensor
category, where $I$ is the identity object and $C$ is the braiding,
its inverse is $C^{-1}$. If $f: U\rightarrow V$, $g: V\rightarrow
W$, $h: I\rightarrow V$, $k: U\rightarrow I$, $\alpha : U \otimes
V\rightarrow P$, $\alpha_{I}: U \otimes V\rightarrow I$ are
morphisms in $\mathcal{D}$, we denote them by:
\[
f=\step
\begin{tangle}
\object{U}\\
\vstr{40}\id\\
\O {f}\\
\vstr{40}\id\\
\object{V}\\
\end{tangle}
\step,\step gf =\step
\begin{tangle}
\object{U}\\
\O {f}\\
\O {g}\\
\object{W}\\
\end{tangle}
\step,\step h=\step
\begin{tangle}
\Q h\\
\object{V}\\
\end{tangle}
\step , \step k=\step
\begin{tangle}
\object{U}\\
\QQ k\\
\end{tangle} \step ,\step
\alpha=\step
\begin{tangle}
\object{U}\step[2]\object{V}\\
\tu \alpha\\
\step\object{P}\\
\end{tangle}
\step,\step \alpha_{I}=\step
\begin{tangle}
\object{U}\step[2]\object{V}\\
\coro {\alpha_{I}}\\
\end{tangle}
\step,\]\[\step C_{U,V}=\step
\begin{tangle}
\object{U}\step[2]\object{V}\\
\x\\
\object{V}\step[2]\object{U}\\
\end{tangle}
\step,\step C_{U,V}^{-1}=\step
\begin{tangle}
\object{V}\step[2]\object{U}\\
\xx\\
\object{U}\step[2]\object{V}\\
\end{tangle}
\step ,\step C_{U,V}=C_{U,V}^{-1}=\step
\begin{tangle}
\object{U}\step[2]\object{V}\\
\XX\\
\object{V}\step[2]\object{U}\\
\end{tangle} \ \ ,\]
where $U,V,W$ are in $\mathcal{D}.$

Since every braided tensor category is always equivalent to a strict
braided tensor category, we can view every braided tensor as a
strict braided tensor and use braiding diagrams freely.

\section {Duality theorem for braided Hopf algebras}\label {e1}

In this section, we obtain the duality theorem for braided  Hopf
algebras living in the braided tensor category $({\cal D}, C).$
Although the most results in  this section appeared in \cite
{ZZH03}, we write them by means of braided diagrams.

If $U,V,W\in ob\mathcal{D}$ and $f,g, act$ are morphisms in
$\mathcal{D}$, we call $act:Hom(V,W)\otimes V\rightarrow W$ satisfy
elimination, if:
\[
\begin{tangle}
\object{U}\step[2]\object{V}\\
\O f\step[2]\id\\
\tu {act}\\
\step\object{W}
\end{tangle}
\step=\step
\begin{tangle}
\object{U}\step[2]\object{V}\\
\O g\step[2]\id\\
\tu {act}\\
\step\object{W}
\end{tangle}\step\Rightarrow f=g.
\]

\begin {Definition} \label {1.1}
Let $(H, m, \eta , \Delta , \epsilon )$ is a braided Hopf algebra in
braided tensor category ${\cal D}$. If there is a braided Hopf
algebra $H^d$ in $\mathcal{D}$ and a morphism $<,>: H^d\otimes
H\rightarrow I$ in $\mathcal{D}$ satisfy:

i)
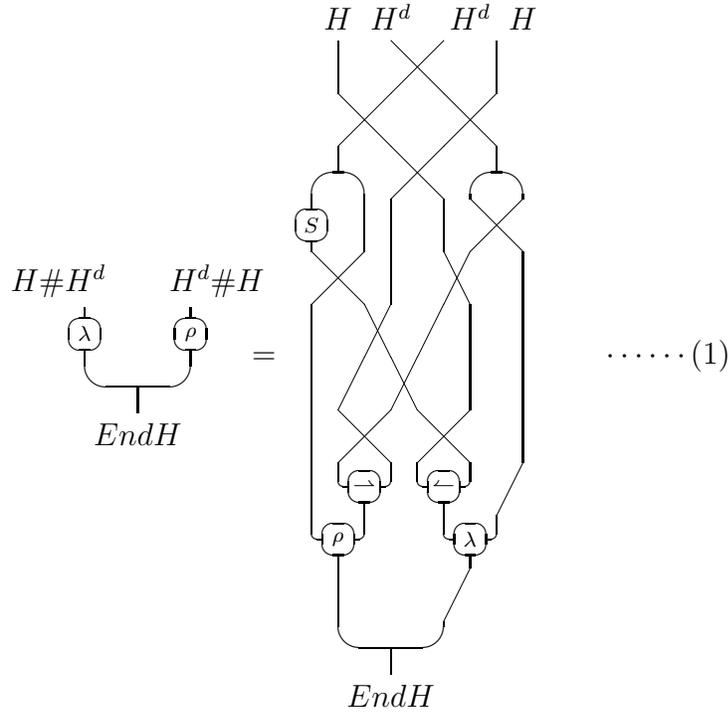
\[
\begin{tangle}
\object{H^{d}}\step\object{H}\step[2]\object{H}\\
\id\step\id\step[2]\id\\
\id\step\cu\\
\coro {<,>}\\
\end{tangle}
\step=\step
\begin{tangle}
\step\object{H^{d}}\step[2]\object{H}\step[2]\object{H}\\
\cd\step\id\step[2]\id\\
\id\step[2]\hx\step[2]\id\\
\coro {<,>}\step\coro{<,>}
\end{tangle}\ \
  ,
\begin{tangle}
\step\object{H^{d}}\\
\step\id\step[2]\Q \eta\\
\step\coro {<,>}
\end{tangle}
\step=\step
\begin{tangle}
\object{H^{d}}\\
\id\\
\QQ \epsilon
\end{tangle}\step;\ \
\]

ii)\[
\begin{tangle}
\object{H^{d}}\step[2]\object{H^{d}}\step[1]\object{H}\\
\id\step[2]\id\step\id\\
\cu \step[1]\id\\
\step \coro {<,>}\\
\end{tangle}
\step=\step
\begin{tangle}
\object{H^{d}}\step[2]\object{H^{d}}\step[2]\object{H}\\
\id\step[2]\id\step\cd\\
\id\step[2]\hx\step[2]\id\\
\coro  {<,>} \step\coro  {<,>}
\end{tangle}\ \
  ,
\step\begin{tangle}
\step[3]\object{H}\\
\step\Q {\eta_{H^{d}}}\step[2]\id\\
\step\coro {<,>}
\end{tangle}
\step=\step
\begin{tangle}
\object{H}\\
\id\\
\QQ \epsilon
\end{tangle}\step; \]

iii)\[
\begin{tangle}
\object{H^{d}}\step[2]\object{H}\\
\O S \step[2]\id\\
\coro {<,>}
\end{tangle}
\step=\step
\begin{tangle}
\object{H^{d}}\step[2]\object{H}\\
\id  \step[2]\O S\\
\coro {<,>}
\end{tangle}
\]
then $H^d$ is called a quasi-dual of $H$ under $<,>$, and $<,>$ is
called a quasi-evaluation of $H^d$ on $H$.
\end {Definition}

No other statement, all the objects and morphisms of this paper are
in $\mathcal{D}$. The $act$ always exists and satisfy elimination,
the braiding  is symmetric on $H$ and $H^d$.
\begin {Lemma} \label {1.3}
(i)  $(H^{d}, \rightharpoonup )$  is a left $H$-module algebra ; \\
(ii)  $(H, \rightharpoonup )$  is a left  $H^{ d}$-module algebra;\\
(iii) $(H^{d}, \leftharpoonup )$  is a right $H$-module algebra;\\
(iv)$(H, \leftharpoonup )$  is a right $H^{d}$-module algebra;\\
where
\[
\begin{tangle}
\object{H}\step[4]\object{H^{d}}\\
\d \step[2]\dd\\
\step\tu \rightharpoonup\\
\step[2]\object{H^{d}}
\end{tangle}
\step=\step
\begin{tangle}
\object{H}\step[3.5]\object{H^{d}}\\
\id \step[2]\cd\\
\XX\step[2]\id\\
\id\step[2]\XX\\
\object{H^{d}}\step[2]\coro {<,>}\\
\end{tangle}\step ;
\step
\begin{tangle}
\object{H^{d}}\step[4]\object{H}\\
\d \step[2]\dd\\
\step\tu \rightharpoonup\\
\step[2]\object{H}
\end{tangle}
\step=\step
\begin{tangle}
\object{H^{d}}\step[3]\object{H}\\
\id \step[2]\cd\\
\XX\step[2]\id\\
\id\step[2]\coro {<,>}\\
\object{H}
\end{tangle}\step;\]\[
\begin{tangle}
\object{H^{d}}\step[4]\object{H}\\
\d \step[2]\dd\\
\step\tu \leftharpoonup\\
\step[2]\object{H^{d}}
\end{tangle}
\step=\step
\begin{tangle}
\step\object{H^{d}}\step[3]\object{H}\\
\cd\step[2]\id\\
\id\step[2]\XX\\
\coro {<,>}\step[2]\id\\
\step[4]\object{H^{d}}
\end{tangle}\step;\step
\begin{tangle}
\object{H}\step[4]\object{H^{d}}\\
\d \step[2]\dd\\
\step\tu \leftharpoonup\\
\step[2]\object{H^{d}}
\end{tangle}
\step=\step
\begin{tangle}
\step\object{H}\step[3]\object{H^{d}}\\
\cd\step[2]\id\\
\id\step[2]\XX\\
\XX\step[2]\id\\
\coro {<,>}\step[2]\object{H}\\
\end{tangle}\step.
\]
\end {Lemma}

Consequently, we construct two smash products $H\#H^d$ and $H^d\#H$
\cite [Chapter 4] {Zh99} .

\begin {Definition} \label  {1.4} We say $CRL$-condition
holds on $H$ and $H^d$ under $<,
>$ if the following conditions are satisfied:\\
$\mathbf{CRL1}$ \ \ \ $E =: End _{\cal C} \ H\in ob\mathcal{D}$ is an algebra under multiplication of composition in $\mathcal{D}$ and satisfy:\\
\[
\begin{tangle}
\object{E}\step[2] \object{E} \step\object{H} \\
\tu {m} \step\id\\
\step\tu {act}\\
\step[2]\object{H}\\
\end{tangle}
\step=\step
\begin{tangle}
\object{E}\step \object{E} \step[2]\object{H}  \\
\id\step\tu {act}\\
\tu {act}\\
\step\object{H}\\
\end{tangle} \step,\step
\begin{tangle}
\step[2]\object{H}\\
\Q {\eta_E}\step[2]\id \\
\tu {act}\\
\step\object{H}\\
\end{tangle}
\step=\step
\begin{tangle}
\object{H}\\
\id \\
\id \\
\object{H}
\end{tangle}\step;
\]
$\mathbf{CRL2}$ \ \ \ There are two morphisms
 $\lambda :  H  \# H ^d \rightarrow E$ and $\rho : H ^d \# H \rightarrow E$
 such that
\[
\begin{tangle}
\object{H\#H^{d}}\step[2]\object{H}\\
\O \lambda \step[2]\id\\
\tu {act} \\
\step\object{H}
\end{tangle}
\step=\step
\begin{tangle}
\object{H}\step[2]\object{H^{d}}\step[3]\object{H}\\
\id\step[2]\id\step[2]\cd\\
\id\step[2]\XX\step[2]\id\\
\cu\step[2]\coro {<,>}\\
\step\object{H}
\end{tangle}
\step \hbox{and}\step
\begin{tangle}
\object{H^{d}\#H}\step[2]\object{H}\\
\O \rho \step[2]\id\\
\tu {act} \\
\step\object{H}
\end{tangle}
\step=\step
\begin{tangle}
\object{H^{d}}\step[2]\object{H}\step[3]\object{H}\\
\id\step[2]\id\step[2]\cd\\
\id\step[2]\XX\step[2]\id\\
\coro {<,>}\step[2]\XX\\
\step[4]\cu\\
\step[5]\object{H}
\end{tangle}\step;
\]
$\mathbf{CRL3}$ \ \  \ there exists a algebra morphism $\bar \lambda
$ from $E$ to $H\# H^d$ such that: $\bar \lambda
\lambda = id _{H \# H^d} $; \\
\end {Definition}
\begin {Lemma}  \label {1.5}
If $CRL$-condition holds on $H$ and $H^d$ under $<,
>$, then $\lambda $ is an algebra morphism from $H\#H ^{d}$ to $E$ and $\rho
$ is an anti-algebra morphism from $H ^d\#H$ to $E$.
\end {Lemma}
{\bf Proof.} We only consider $\lambda$. It is straightforward to
prove that $\lambda\eta_{H\#H^{d}}=\eta_{E}$ and $ act ((\lambda
\otimes id_{H}) (m \otimes id_{H}))= act (( m \otimes id_{H})
(\lambda \otimes \lambda \otimes id_{H}))$. $\Box$
\begin {Lemma} \label {1.6}
If $CRL$-condition holds on $H$ and $H^d$ under $<,
>$, then $\lambda $ and $\rho $ satisfy the following:
\[
\begin{tangle}
\step\object{H\# {{H^{d}}}}\step[6]\object{{H^{d}}\#
{H}}\\
\Step\O \lambda\step[4]\O \rho\\
\Step\Cu\\
\step[4]\object{End H}\\
\end{tangle}
\step=\step
\begin{tangle}
\step\object{H}\step[2]\object{H^{d}}\step[3]\object{H^{d}}\Step\object{H}\\
\step\id\Step\XX\Step\id\\
\step\XX\Step\XX\\
\cd\step\XX\step\cd\\
\O S\Step\id\step\id\Step\id\step\XX\\
\XX\step\id\Step\X\Step\id\\
\id\Step\X\step\dd\step\id\Step\id\\
\id\step\dd\step\X\Step\id\Step\id\\
\id\step\XX\step\XX\Step\id\\
\id\step\tu \rightharpoonup\step\tu \leftharpoonup\step\dd\\
\tu \rho\step[3]\tu \lambda\\
\step\id\step[4]\dd\\
\step\Cu\\
\step[3]\object{End H}\\
\end{tangle}
\ \ \    \ \ \ \ \cdots\cdots(1)
\]
\end {Lemma}

{\bf Proof.} We show (1) by following five steps. It is easy to
check the following (i) and (ii).

(i)
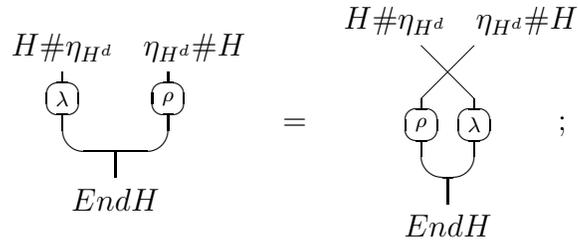
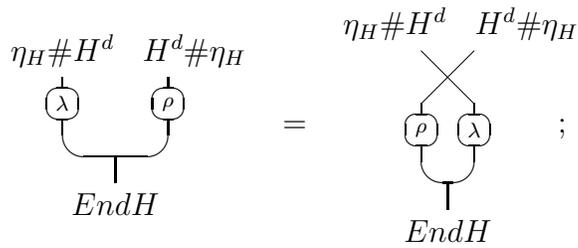
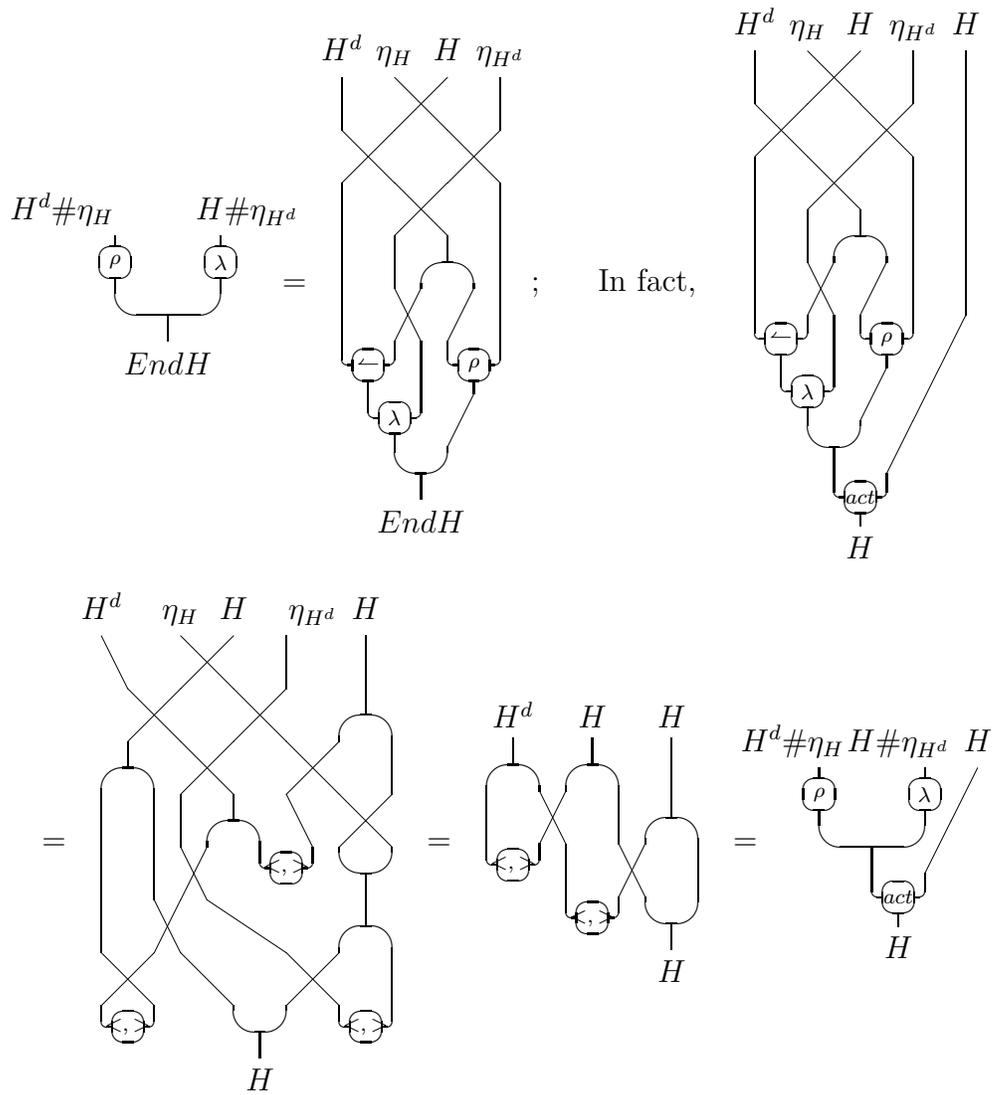
\[
\begin{tangle}
\step[2]\object{H\# {\eta_{H^{d}}}}\step[5]\object{{\eta_{H^{d}}}\#
H}\\
\Step\O \lambda\step[4]\O \rho\\
\Step\Cu\\
\step[4]\object{End H}\\
\end{tangle}
\step[3]=\step[3]
\begin{tangle}
\object{H\# {\eta_{H^{d}}}}\step[5]\object{{\eta_{H^{d}}}\#
H}\\
\step\XX\\
\step\O \rho\step[2]\O \lambda\\
\step\cu\\
\step[2]\object{End H}\\
\end{tangle}\ \ \ ;
\]

(ii)
\[
\begin{tangle}
\step[2]\object{\eta_H\# {H^{d}}}\step[5]\object{H^{d}\#
\eta_H}\\
\Step\O \lambda\step[4]\O \rho\\
\Step\Cu\\
\step[4]\object{End H}\\
\end{tangle}
\step[3]=\step[3]
\begin{tangle}
\object{\eta_H\# {H^{d}}}\step[5]\object{H^{d}\#
\eta_H}\\
\step\XX\\
\step\O \rho\step[2]\O \lambda\\
\step\cu\\
\step[2]\object{End H}\\
\end{tangle}\ \ \ ;
\]

(iii)
\[
\begin{tangle}
\object{{H^{d}}\#\eta_H}\step[7]\object{{H}\#\eta_{H^{d}}
}\\
\Step\O \rho\step[4]\O \lambda\\
\Step\Cu\\
\step[4]\object{End H}\\
\end{tangle}
\step=\step
\begin{tangle}
\object{H^{d}}\step[2]\object{\eta_H}\step[2]\object{H}\step[2]\object{\eta_{H^{d}}}\\
\id\Step\XX\Step\id\\
\XX\Step\XX\\
\id\Step\XX\Step\id\\
\id\Step\id\step\cd\step\id\\
\id\Step\X\step\dd\step\id\\
\tu \leftharpoonup\step\id\step\tu \rho\\
\step\tu \lambda\step\dd\\
\Step\cu\\
\step[3]\object{End H}\\
\end{tangle}\ \ \ ;
\step[2]\hbox{In fact},\step[2]
\begin{tangle}
\object{H^{d}}\step[2]\object{\eta_H}\step[2]\object{H}\step[2]\object{\eta_{H^{d}}}\step[2]\object{H}\\
\id\Step\XX\Step\id\Step\id\\
\XX\Step\XX\Step\id\\
\id\Step\XX\Step\id\Step\id\\
\id\Step\id\step\cd\step\id\Step\id\\
\id\Step\X\step\dd\step\id\Step\id\\
\tu \leftharpoonup\step\id\step\tu \rho\step\dd\\
\step\tu \lambda\step\dd\step\dd\\
\Step\cu\step\dd\\
\step[3]\tu {act}\\
\step[4]\object{ H}\\
\end{tangle} \]
\[\step=\step
\begin{tangle}
\object{H^{d}}\step[3]\object{\eta_H}\step[2]\object{H}\step[3]\object{\eta_{H^{d}}}\step[2]\object{H}\\
\d\Step\XX\Step\id\step[3]\id\\
\step\XX\Step\XX\Step\cd\\
\cd\step\XX\Step\XX\Step\id\\
\id\Step\id\step\id\step\cd\step\d\step\XX\\
\id\Step\id\step\X\Step\coro {<,>}\step\cu\\
\id\Step\X\step\nw3\step[4]\cd\\
\XX\step\nw2\step[3]\XX\Step\id\\
\coro {<,>}\step[3]\cu\Step\coro {<,>}\\
\step[6]\object{H}\\
\end{tangle}
\step=\step
\begin{tangle}
\step\object{H^{d}}\step[3]\object{H}\step[3]\object{H}\\
\cd\step\cd\Step\id\\
\id\Step\X\Step\id\step\cd\\
\coro {<,>}\step\id\Step\X\Step\id\\
\step[3]\coro {<,>}\step\cu\\
\step[7]\object{H}\\
\end{tangle}
\step=\step \begin{tangle}
\object{{H^{d}}\#\eta_H}\step[4]\object{{H}\#\eta_{H^{d}}
}\step[3]\object{H}\\
\step\O \rho\step[4]\O \lambda\step\dd\\
\step\Cu\dd\\
\step[3]\tu {act}\\
\step[4]\object{ H}\\
\end{tangle}\step.
\]
Thus (iii) holds.

(iv)
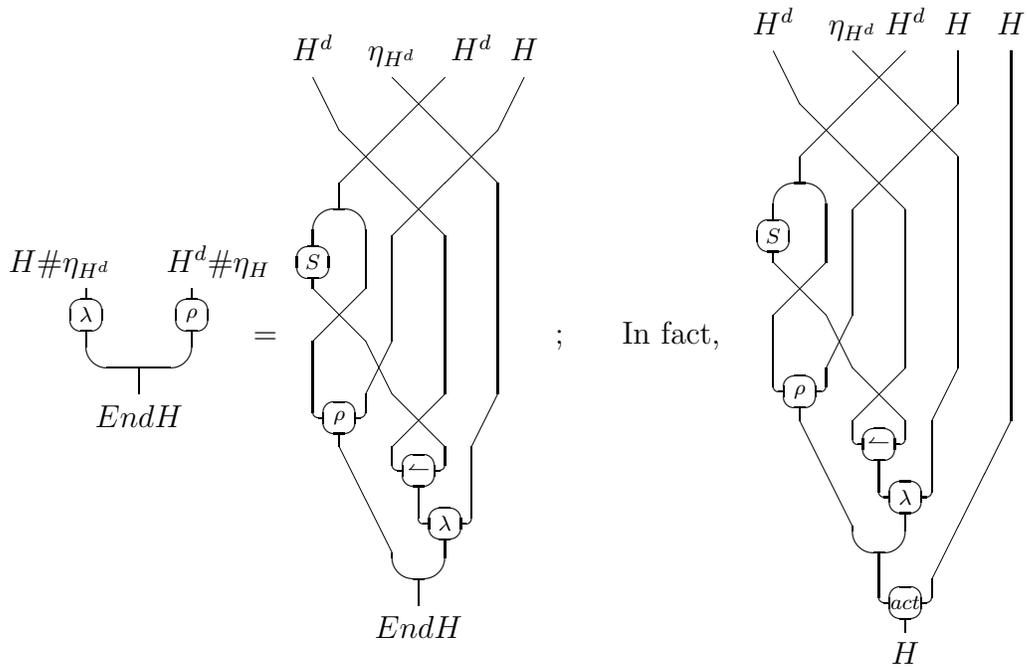
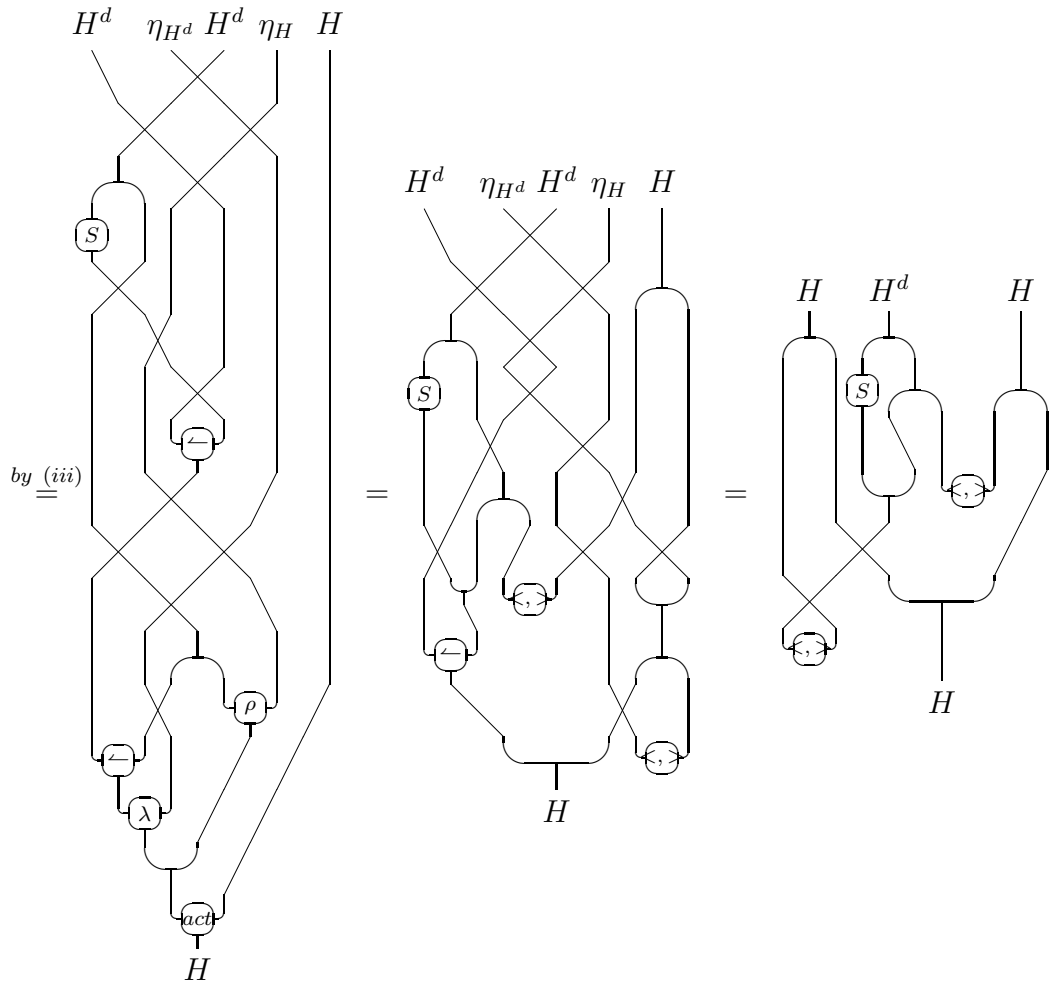
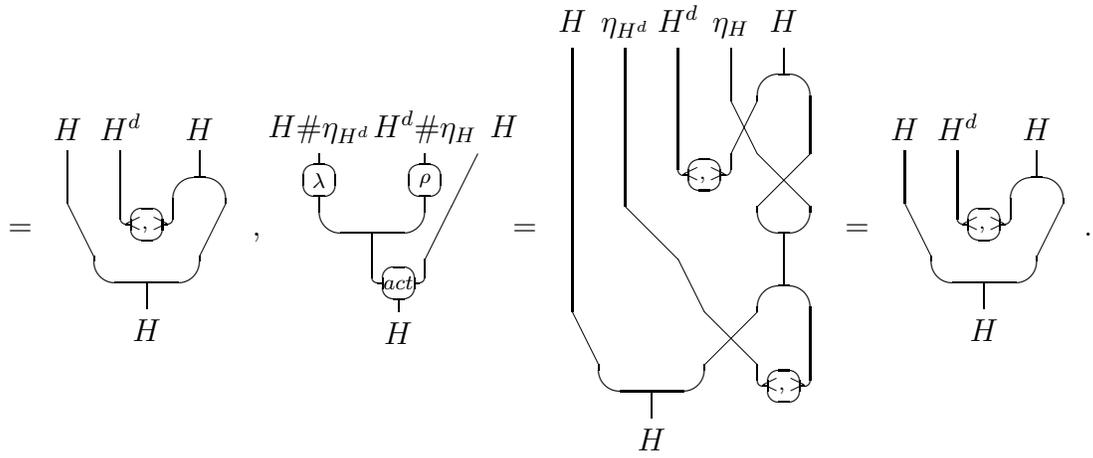
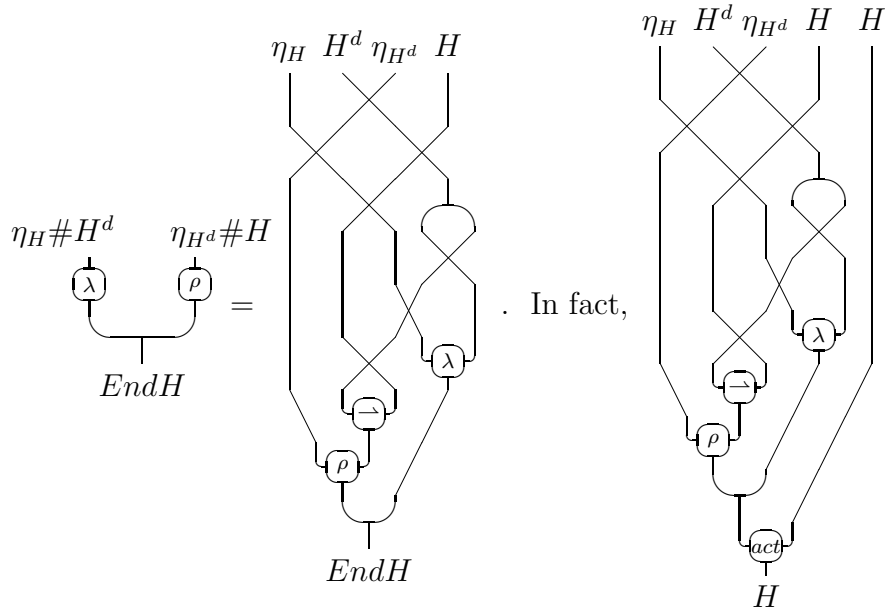
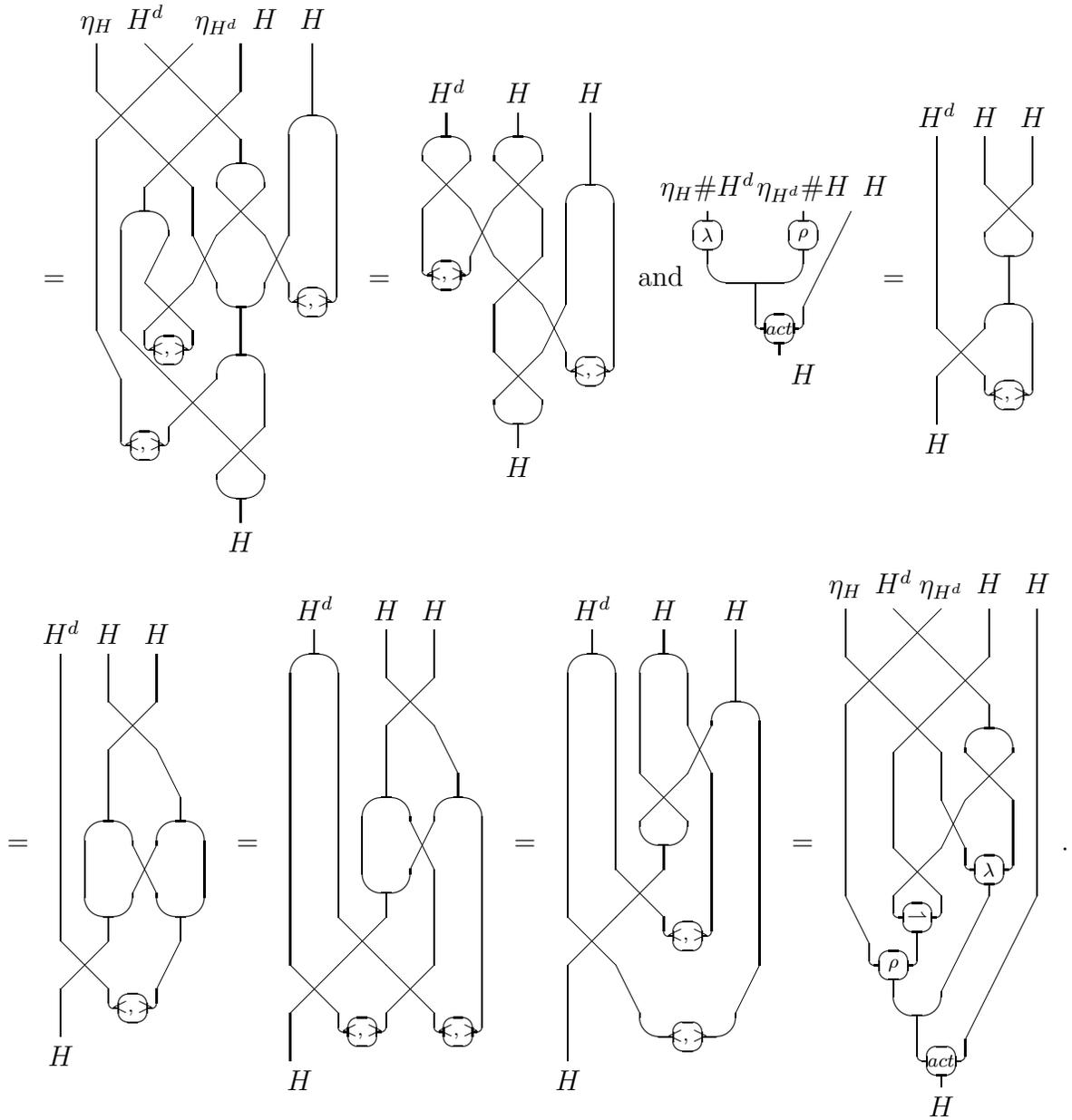
\[
\begin{tangle}
\step\object{H\# {\eta_{H^{d}}}}\step[6]\object{{H^{d}}\#
{\eta_H}}\\
\Step\O \lambda\step[4]\O \rho\\
\Step\Cu\\
\step[4]\object{EndH}\\
\end{tangle}
\step=\step
\begin{tangle}
\object{H^{d}}\step[3]\object{\eta_{H^{d}}}\step[3]\object{H^{d}}\step[2]\object{H}\\
\d\Step\XX\Step\dd\\
\step\XX\Step\XX\\
\cd\step\XX\Step\id\\
\O S\Step\id\step\id\Step\id\Step\id\\
\XX\step\id\Step\id\Step\id\\
\id\Step\X\Step\id\Step\id\\
\tu \rho\step\XX\step\dd\\
\step\d\step\tu \leftharpoonup\step\id\\
\Step\d\step\tu \lambda\\
\step[3]\cu\\
\step[4]\object{End H}\\
\end{tangle}\ \ \ ;
\step[2] \hbox{In fact,}\step[2]
\begin{tangle}
\object{H^{d}}\step[3]\object{\eta_{H^{d}}}\step[2]\object{H^{d}}\step[2]\object{H}\step[2]\object{H}\\
\d\Step\XX\Step\id\Step\id\\
\step\XX\Step\XX\Step\id\\
\cd\step\XX\Step\id\Step\id\\
\O S\Step\id\step\id\Step\id\Step\id\Step\id\\
\XX\step\id\Step\id\Step\id\Step\id\\
\id\Step\X\Step\id\Step\id\Step\id\\
\tu \rho\step\XX\step\dd\Step\id\\
\step\d\step\tu \leftharpoonup\step\id\step[2]\dd\\
\Step\d\step\tu \lambda\step\dd\\
\step[3]\cu\step\dd\\
\step[4]\tu {act}\\
\step[5]\object{H}
\end{tangle}\]
\[\stackrel{by\ (iii)}{=}
\begin{tangle}
\object{H^{d}}\step[3]\object{\eta_{H^{d}}}\step[2]\object{H^{d}}\step[2]\object{\eta_{H}}\step[2]\object{H}\\
\d\Step\XX\Step\id\step[2]\id\\
\step\XX\Step\XX\step[2]\id\\
\cd\step\XX\Step\id\step[2]\id\\
\O S\Step\id\step\id\Step\id\Step\id\step[2]\id\\
\XX\step\id\Step\id\Step\id\step[2]\id\\
\id\Step\X\Step\id\Step\id\step[2]\id\\
\id\Step\id\step\XX\Step\id\step[2]\id\\
\id\Step\id\step\tu \leftharpoonup\Step\id\step[2]\id\\
\id\Step\XX\Step\dd\step[2]\id\\
\XX\Step\XX\step[3]\id\\
\id\Step\XX\Step\d\step[2]\id\\
\id\Step\id\step\cd\Step\id\step[2]\id\\
\id\Step\X\Step\tu \rho\step\dd\\
\tu \leftharpoonup\step\id\Step\dd\step\dd\\
\step\tu \lambda\step\dd\step\dd\\
\Step\cu\step\dd\\
\step[3]\tu {act}\\
\step[4]\object{ H}\\
\end{tangle}
\step=\step
\begin{tangle}
\object{H^{d}}\step[3]\object{\eta_{H^{d}}}\step[2]\object{H^{d}}\step[2]\object{\eta_{H}}\step[2]\object{H}\\
\d\Step\XX\Step\id\Step\id\\
\step\XX\Step\XX\step\cd\\
\cd\step\XX\Step\id\step\id\Step\id\\
\O S\step[2]\id\step\XX\Step\id\step\id\Step\id\\
\id\step[2]\X\Step\XX\step\id\Step\id\\
\id\step\dd\cd\step\id\Step\X\Step\id\\
\X\step\id\step\dd\step\XX\step\XX\\
\id\step\hdcu\step\coro {<,>}\Step\id\step\cu\\
\tu \leftharpoonup\step[5]\id\step\cd\\
\step\nw2\step[5]\X\Step\id\\
\step[3]\Cu\step\coro {<,>}\\
\step[5]\object{H}\\
\end{tangle}
\step=\step
\begin{tangle}
\step\object{H}\step[3]\object{H^{d}}\step[5]\object{H}\\
\cd\step\cd\step[4]\id\\
\id\Step\id\step\O S\step\cd\Step\cd\\
\id\Step\id\step\id\step\d\step\id\Step\id\Step\id\\
\id\Step\id\step\cu\step\coro {<,>}\step\dd\\
\id\Step\XX\step[4]\dd\\
\XX\Step\Cu\\
\coro {<,>}\step[4]\id\\
\step[6]\object{H}\\
\end{tangle}\]
\[\step=\step
\begin{tangle}
\object{H}\step[2]\object{H^{d}}\step[3]\object{H}\\
\id\Step\id\Step\cd\\
\d\step\coro {<,>}\step\dd\\
\step\Cu\\
\step[3]\object{H}\\
\end{tangle}
\step ,\step[2]
\begin{tangle}
\object{H\# {\eta_{H^{d}}}}\step[4]\object{{H^{d}}\#
{\eta_H}}\step[3]\object{H}\\
\O \lambda\step[4]\O \rho\step\dd\\
\Cu\dd\\
\step[2]\tu {act}\\
\step[3]\object{H}\\
\end{tangle}=\step
\begin{tangle}
\object{H}\step[2]\object{\eta_{H^{d}}}\step[2]\object{H^{d}}\step[2]\object{\eta_H}\step[2]\object{H}\\
\id\Step\id\step[2]\id\Step\id\step\cd\\
\id\Step\id\step[2]\id\Step\X\Step\id\\
\id\Step\id\step[2]\coro {<,>}\step\XX\\
\id\step[2]\nw2\step[4]\cu\\
\id\step[4]\d\step[2]\cd\\
\d\step[4]\XX\Step\id\\
\step\Cu\step[2]\coro {<,>}\\
\step[3]\object{H}\\
\end{tangle}
\step=\step
\begin{tangle}
\object{H}\step[2]\object{H^{d}}\step[3]\object{H}\\
\id\Step\id\Step\cd\\
\d\step\coro {<,>}\step\dd\\
\step\Cu\\
\step[3]\object{H}\\
\end{tangle}\ \ .
\]

Thus (iv) holds.

(v)
\[
\begin{tangle}
\step\object{\eta_H\# {H^{d}}}\step[6]\object{\eta_{H^{d}}\#
H}\\
\Step\O \lambda\step[4]\O \rho\\
\Step\Cu\\
\step[4]\object{End H}\\
\end{tangle}
=\step
\begin{tangle}
\object{\eta_H}\Step\object{H^{d}}\step[2]\object{\eta_{H^{d}}}\step[2]\object{H}\\
\id\Step\XX\Step\id\\
\XX\Step\XX\\
\id\Step\XX\step\cd\\
\id\Step\id\Step\id\step\XX\\
\id\Step\id\Step\X\Step\id\\
\id\Step\XX\step\tu \lambda\\
\d\step\tu \rightharpoonup\step\dd\\
\step\tu \rho\step\dd\\
\Step\cu\\
\step[3]\object{End H}\\
\end{tangle}\step.\ \
\hbox{In fact},\step\begin{tangle}
\object{\eta_H}\Step\object{H^{d}}\step[2]\object{\eta_{H^{d}}}\step[2]\object{H}\step[2]\object{H}\\
\id\Step\XX\Step\id\Step\id\\
\XX\Step\XX\Step\id\\
\id\Step\XX\step\cd\step\id\\
\id\Step\id\Step\id\step\XX\step\id\\
\id\Step\id\Step\X\Step\id\step\id\\
\id\Step\XX\step\tu \lambda\step\id\\
\d\step\tu \rightharpoonup\step\dd\step\dd\\
\step\tu \rho\step\dd\step\dd\\
\Step\cu\step\dd\\
\step[3]\tu {act}\\
\step[4]\object{H}
\end{tangle}\]
\[
=\step
\begin{tangle}
\object{\eta_H}\Step\object{H^{d}}\step[3]\object{\eta_{H^{d}}}\step[2]\object{H}\step[2]\object{H}\\
\id\Step\XX\Step\id\step[3]\id\\
\XX\Step\XX\Step\cd\\
\id\Step\XX\step\cd\step\id\Step\id\\
\id\step\cd\step\id\step\XX\step\id\Step\id\\
\id\step\id\step\dd\step\X\Step\X\Step\id\\
\id\step\id\step\XX\step\cu\step\coro {<,>}\\
\d\nw2\coro {<,>}\step\cd\\
\step\id\Step\XX\Step\id\\
\step\coro {<,>}\Step\XX\\
\step[5]\cu\\
\step[6]\object{H}\\
\end{tangle}
\step=\step
\begin{tangle}
\step\object{H^{d}}\step[3]\object{H}\step[3]\object{H}\\
\cd\step\cd\Step\id\\
\XX\step\XX\step\cd\\
\id\Step\X\Step\id\step\id\Step\id\\
\coro {<,>}\step\XX\step\id\Step\id\\
\step[3]\id\Step\X\Step\id\\
\step[3]\XX\step\coro {<,>}\\
\step[3]\cu\\
\step[4]\object{H}\\
\end{tangle}
\step \hbox{and}\step
\begin{tangle}
\object{\eta_H\# {H^{d}}}\step[4]\object{\eta_{H^{d}}\#
H}\step[3]\object{H}\\
\O \lambda\step[4]\O \rho\step\dd\\
\Cu\dd\\
\step[2]\tu {act}\\
\step[4]\object{ H}\\
\end{tangle}=\step
\begin{tangle}
\object{H^{d}}\step[2]\object{H}\step[2]\object{H}\\
\id\Step\id\Step\id\\
\id\Step\XX\\
\id\Step\cu\\
\id\Step\cd\\
\XX\Step\id\\
\id\Step\coro {<,>}\\
\object{H}\\
\end{tangle}\]
\[=\step
\begin{tangle}
\object{H^{d}}\step[2]\object{H}\step[2]\object{H}\\
\id\Step\id\Step\id\\
\id\Step\XX\\
\id\Step\id\Step\d\\
\id\step\cd\step\cd\\
\id\step\id\Step\X\Step\id\\
\id\step\cu\step\cu\\
\XX\Step\dd\\
\id\Step\coro {<,>}\\
\object{H}\\
\end{tangle}
\step=\step
\begin{tangle}
\step\object{H^{d}}\step[3]\object{H}\step[2]\object{H}\\
\cd\step[2]\id\Step\id\\
\id\Step\id\step[2]\XX\\
\id\step[2]\id\Step\id\Step\d\\
\id\step[2]\id\step\cd\step\cd\\
\id\Step\id\step\id\Step\X\Step\id\\
\id\Step\id\step\cu\step\id\Step\id\\
\id\Step\XX\Step\id\Step\id\\
\XX\Step\XX\Step\id\\
\id\Step\coro {<,>}\Step\coro {<,>}\\\
\object{H}\\
\end{tangle}
\step=\step
\begin{tangle}
\step\object{H^{d}}\step[3]\object{H}\step[3]\object{H}\\
\cd\step\cd\Step\id\\
\id\Step\id\step\id\Step\id\step\cd\\
\id\Step\id\step\id\Step\X\Step\id\\
\id\Step\id\step\XX\step\id\Step\id\\
\id\Step\id\step\cu\step\id\Step\id\\
\id\Step\XX\Step\id\step[2]\id\\
\XX\Step\coro {<,>}\step[2]\id\\
\id\Step\d\step[4]\dd\\
\id\step[3]\coRo {<,>}\\
\object{H}\\
\end{tangle}
\step=\step\begin{tangle}
\object{\eta_H}\Step\object{H^{d}}\step[2]\object{\eta_{H^{d}}}\step[2]\object{H}\step[2]\object{H}\\
\id\Step\XX\Step\id\Step\id\\
\XX\Step\XX\Step\id\\
\id\Step\XX\step\cd\step\id\\
\id\Step\id\Step\id\step\XX\step\id\\
\id\Step\id\Step\X\Step\id\step\id\\
\id\Step\XX\step\tu \lambda\step\id\\
\d\step\tu \rightharpoonup\step\dd\step\dd\\
\step\tu \rho\step\dd\step\dd\\
\Step\cu\step\dd\\
\step[3]\tu {act}\\
\step[4]\object{ H}\\
\end{tangle}\step.
\]
Thus (v) holds. Next, we show (1) holds.
\[
\hbox{the left side of (1)}\stackrel{\hbox{by Lemma \ref {1.5}
}}{\step=\step[5]}
\begin{tangle}
\object{H\# {\eta_{H^{d}}}}\step[5]\object{{\eta_H}\#
H^{d}}\step[5]\object{H^{d}\#\eta_H}\step[5]\object{\eta_{H^{d}}\# H}\\
\O \lambda\step[3]\Step\O \lambda\step[6]\XX \\
\id\step[5]\id\step[6]\O \rho\step[2]\O \rho\\
\nw3\step[4]\nw3\step[5]\cu\\
\step[3]\nw3\step[4]\Cu\\
\step[6]\Cu\\
\step[8]\object{End H}\\
\end{tangle}
\step \]\[\stackrel{\hbox{by (v)}}{=}\step
\begin{tangle}
\object{H\#
{\eta_{H^{d}}}}\step[4]\object{\eta_H}\step[3]\object{H^{d}}\step[3]
\object{H^{d}\#\eta_H}\step[4]\object{\eta_{H^{d}}}\step[2]\object{H}\\
\O \lambda\step[4]\nw2\step[2]\nw2\step[3]\XX\step[2]\dd\\
\d\step[5]\d\Step\XX\step[2]\XX\\
\step\d\step[5]\XX\Step\XX\step[2]\O \rho\\
\step[2]\d\step[4]\id\Step\XX\step\cd\step\id\\
\step[3]\d\step[3]\id\Step\id\Step\id\step\XX\step\id\\
\step[4]\d\step[2]\id\Step\id\Step\X\Step\id\step\id\\
\step[5]\d\step\id\Step\XX\step\tu \lambda\step\id\\
\step[6]\id\step\d\step\tu \rightharpoonup\step[2]\cu\\
\step[6]\id\step[2]\tu \rho\step[3]\dd\\
\step[6]\d\step[2]\Cu\\
\step[7]\Cu\\
\step[9]\object{End H}\\
\end{tangle}\]\[
\step\stackrel{\hbox{by (i)(ii)}}{=}\step
\begin{tangle}
\object{H}\step[2]\object{\eta_{H^{d}}}\step[2]\object{\eta_H}\step[2]\object{H^{d}}
\step[2]\object{H^{d}}\step[2]\object{\eta_H}\step[2]\object{\eta_{H^{d}}}\step[2]\object{H}\\
\id\Step\id\Step\id\Step\id\Step\id\Step\XX\Step\id\\
\id\Step\id\Step\id\Step\id\Step\XX\Step\XX\\
\id\Step\id\Step\id\Step\XX\Step\XX\Step\id\\
\id\Step\id\Step\XX\Step\XX\Step\id\Step\id\\
\id\Step\XX\Step\XX\step\cd\step\id\Step\id\\
\XX\Step\id\Step\id\Step\id\step\XX\step\id\Step\id\\
\id\Step\id\Step\id\Step\id\Step\X\Step\X\Step\id\\
\id\Step\id\Step\id\Step\XX\step\XX\step\XX\\
\id\Step\id\Step\d\step\tu
\rightharpoonup\step\id\Step\X\Step\id\\
\id\Step\d\Step\XX\Step\tu \rho\step\tu \lambda\\
\d\Step\XX\Step\id\step[3]\id\Step\dd\\
\step\tu \rho\Step\tu \lambda\step[3]\id\Step\id\\
\step[2]\nw2\step[3]\Cu\Step\id\\
\step[4]\nw2\step[3]\Cu\\
\step[6]\Cu\\
\step[8]\object{End H}\\
\end{tangle}
\step\stackrel{\hbox{by (iv)}}{=}\step
\begin{tangle}
\object{H}\step[2]\object{\eta_{H^{d}}}\step[2]\object{\eta_H}\step[2]\object{H^{d}}
\step[2]\object{H^{d}}\step[2]\object{\eta_H}\step[2]\object{\eta_{H^{d}}}\step[2]\object{H}\\
\id\Step\id\Step\id\Step\id\Step\id\Step\XX\Step\id\\
\id\Step\id\Step\id\Step\id\Step\XX\Step\XX\\
\id\Step\id\Step\id\Step\XX\Step\XX\Step\id\\
\id\Step\id\Step\XX\Step\XX\Step\id\Step\id\\
\id\Step\XX\Step\XX\step\cd\step\id\Step\id\\
\XX\Step\id\Step\id\Step\id\step\XX\step\id\Step\id\\
\id\Step\id\Step\id\Step\id\Step\X\Step\X\Step\id\\
\id\Step\id\Step\id\Step\XX\step\XX\step\XX\\
\id\Step\id\Step\d\step\tu
\rightharpoonup\step\id\Step\X\Step\id\\
\id\Step\d\Step\XX\Step\id\Step\id\step\tu \lambda\\
\d\Step\XX\Step\XX\Step\id\Step\id\\
\step\tu \rho\Step\XX\Step\XX\Step\id\\
\step[2]\id\Step\cd\step\XX\Step\id\Step\id\\
\step[2]\id\Step\O S\Step\id\step\id\Step\id\Step\id\Step\id\\
\step[2]\id\Step\XX\step\id\Step\id\Step\id\Step\id\\
\step[2]\id\Step\id\Step\X\Step\id\Step\id\Step\id\\
\step[2]\d\step\tu \rho\step\XX\Step\id\Step\id\\
\step[3]\cu\Step\tu \leftharpoonup\step\dd\Step\id\\
\step[4]\d\step[3]\tu \lambda\step[2]\dd\\
\step[5]\Cu\step[2]\dd\\
\step[7]\Cu\\
\step[9]\object{End H}\\
\end{tangle}\]\[
\step\stackrel{\hbox{by Lemma \ref {1.5}  }}{=}
\begin{tangle}
\object{H}\step[3]\object{H^{d}}\step[2]\object{H^{d}}\Step\object{H}\\
\id\step[3]\id\Step\XX\\
\id\step[3]\XX\Step\id\\
\id\step[3]\id\step\cd\step\id\\
\id\Step\dd\step\XX\step\id\\
\id\Step\XX\Step\X\\
\d\step\tu \rightharpoonup\Step\id\step\id\\
\step\XX\Step\dd\step\id\\
\dd\Step\XX\Step\id\\
\id\Step\cd\step\d\step\id\\
\id\Step\O S\Step\id\Step\id\step\id\\
\id\Step\XX\Step\id\step\id\\
\XX\Step\XX\step\id\\
\tu \rho\Step\tu \leftharpoonup\step\id\\
\step\d\step[3]\tu \lambda\\
\Step\Cu\\
\step[4]\object{End H}\\
\end{tangle}
\step=\step
\begin{tangle}
\step\object{H}\step[2]\object{H^{d}}\step[3]\object{H^{d}}\Step\object{H}\\
\step\id\Step\XX\Step\id\\
\step\XX\Step\XX\\
\cd\step\XX\step\cd\\
\O S\Step\id\step\id\Step\id\step\XX\\
\XX\step\id\Step\X\Step\id\\
\id\Step\X\step\dd\step\id\Step\id\\
\id\step\dd\step\X\Step\id\Step\id\\
\id\step\XX\step\XX\Step\id\\
\id\step\tu \rightharpoonup\step\tu \leftharpoonup\step\dd\\
\tu \rho\step[3]\tu \lambda\\
\step\id\step[4]\dd\\
\step\Cu\\
\step[3]\object{End H}\\
\end{tangle}
\step=\step
\begin{tangle}
\hbox{the right side of (1). }
\end{tangle}\Box
\]

If $(R , \psi )$ is a right $H^{d}$-comodule algebra , then $(R,
\alpha )$ becomes a left $H$-module algebra (similar to  \cite
[Example 1.6.4] {Ma95b}) under the module operation:
\[
\begin{tangle}
\object{H}\step[2]\object{R}\\
\tu \alpha\\
\step\object{R}\\
\end{tangle}
\step=\step
\begin{tangle}
\object{H}\step[3]\object{R}\\
\id\step[2]\td \psi\\
\x\step[2]\id\\
\id\step[2]\XX\\
\object{R}\step[2]\coro {<>}\\
\end{tangle}\step.
\]
\begin {Lemma} \label {1.7}
$R \# H$ becomes an $H^{ d}$-module algebra under the
 module operation:
\[
\begin{tangle}
\object{H^{d}}\step[4]\object{R\#H}\\
\d \step[2]\dd\\
\step\tu {\rightharpoonup'} \\
\step[2]\object{H^{d}}
\end{tangle}
\step=\step
\begin{tangle}
\object{H^{d}}\step[2]\object{R}\step[2]\object{H}\\
\x\step[2]\id\\
\id\step[2]\tu \rightharpoonup\\
\object{R}\step[3]\object{H}\\
\end{tangle}\step.\]
\end {Lemma}

{\bf Proof.} It is straightforward. $\Box$

Consequently, we obtain another  smash product $(R \# H)\# H^d.$

\begin {Theorem} \label {1.8}
 Let  $H$ and $H^{d}$ be Hopf algebra with invertible antipodes, and the $CRL$-condition holds on $H$ and $H^d$ under
 $<,>$.
Let $R$ be an $H^{d}$-comodule algebra such that  $R$ is an
$H$-module algebra defined  as  above, $H^{d}$ act on $R\#H$ by
acting trivially on $R$ and via $\rightharpoonup $ on $H$, then
  $$   (R \# H)\# H^{d}   \cong R \otimes (H
  \# H^{d} ) \hbox { \ \ \  as }  \hbox { algebras in }  {\cal D}.$$
\end {Theorem}
{\bf Proof. } By (CRL)-condition,  there exists a algebra morphism
$\bar \lambda $ from $E$ to $H\# H^{d}$ such that $\bar \lambda
\lambda = id _{H \# H^{d}}$, We first define a morphism $w=\bar
\lambda \rho (S^{-1} \otimes \eta _H )$ from $H^{d} $ to $H \#
H^{d}$. Since $\rho $ and $S^{-1}$ are anti-algebra morphisms by
Lemma \ref {1.5}, $w$ is an algebra morphism.

We now define two morphisms:
\[
\begin{tangle}
\object{R\otimes(H\# H^{d})}\\
\id\\
\O \Psi\\
\id\\
\object{(R\# H)\# H^{d}}\\
\end{tangle}
\step[3]=\step[3]
\begin{tangle}
\step\object{R}\step[4]\object{H}\Step\object{H^{d}}\\
\td \psi\step[3]\id\Step\id\\
\id\Step\O S\step[3]\id\Step\id\\
\id\step\td w\Step\id\Step\id\\
\id\step\id\step\cd\step\id\Step\id\\
\id\step\id\step\id\Step\X\Step\id\\
\id\step\id\step\tu \rightharpoonup\step\cu\\
\id\step\cu\step[3]\id\\
\object{R}\Step\object{H}\step[4]\object{H^{d}}\\
\end{tangle}
\step[3]and\step[5]
\begin{tangle}
\object{(R\# H)\# H^{d}}\\
\id\\
\O \Phi\\
\id\\
\object{R\otimes(H\# H^{d})}\\
\end{tangle}
\step[3]=\step[3]
\begin{tangle}
\step\object{R}\step[4]\object{H}\Step\object{H^{d}}\\
\td \psi\step[3]\id\Step\id\\
\id\step\td w\Step\id\Step\id\\
\id\step\id\step\cd\step\id\Step\id\\
\id\step\id\step\id\Step\X\Step\id\\
\id\step\id\step\tu \rightharpoonup\step\cu\\
\id\step\cu\step[3]\id\\
\object{R}\Step\object{H}\step[4]\object{H^{d}}\\
\end{tangle}\ \ \ .
\]
It is straightforward to check that $\Psi\Phi=id,\Phi \Psi=id$. Now
we show that $\Phi$ is algebraic. Define
\[
\begin{tangle}
\step\object{(R\# H)\# H^{d}}\\
\step\id\\
\obox 2{\Phi'}\\
\step\id\\
\step\object{End H}\\
\end{tangle}
\step=\step[3]
\begin{tangle}
\step\object{R}\step[5]\object{H\# H^{d}}\\
\td \xi\step[4]\O \lambda\\
\id\Step\Cu\\
\object{R}\step[4]\object{End H}\\
\end{tangle}\step[2]
where \step[2]
\begin{tangle}
\object{R}\\
\id\\
\O \xi\\
\id\\
\object{R\otimes End H}\\
\end{tangle}
\step[3]=\step[3]
\begin{tangle}
\step\object{R}\step[3]\object{\eta_H}\\
\td \psi\step[2]\id\\
\id\step[2]\O {\bar{S}}\step[2]\id\\
\id\Step\tu \rho\\
\object{R}\step[4]\object{End H}\\
\end{tangle}\ \ \ .
\]
It is clear that $\Phi = (id \otimes \bar{\lambda} )\Phi '$.
Consequently, we only need to show that $\Phi '$ is alegebraic. We
claim that

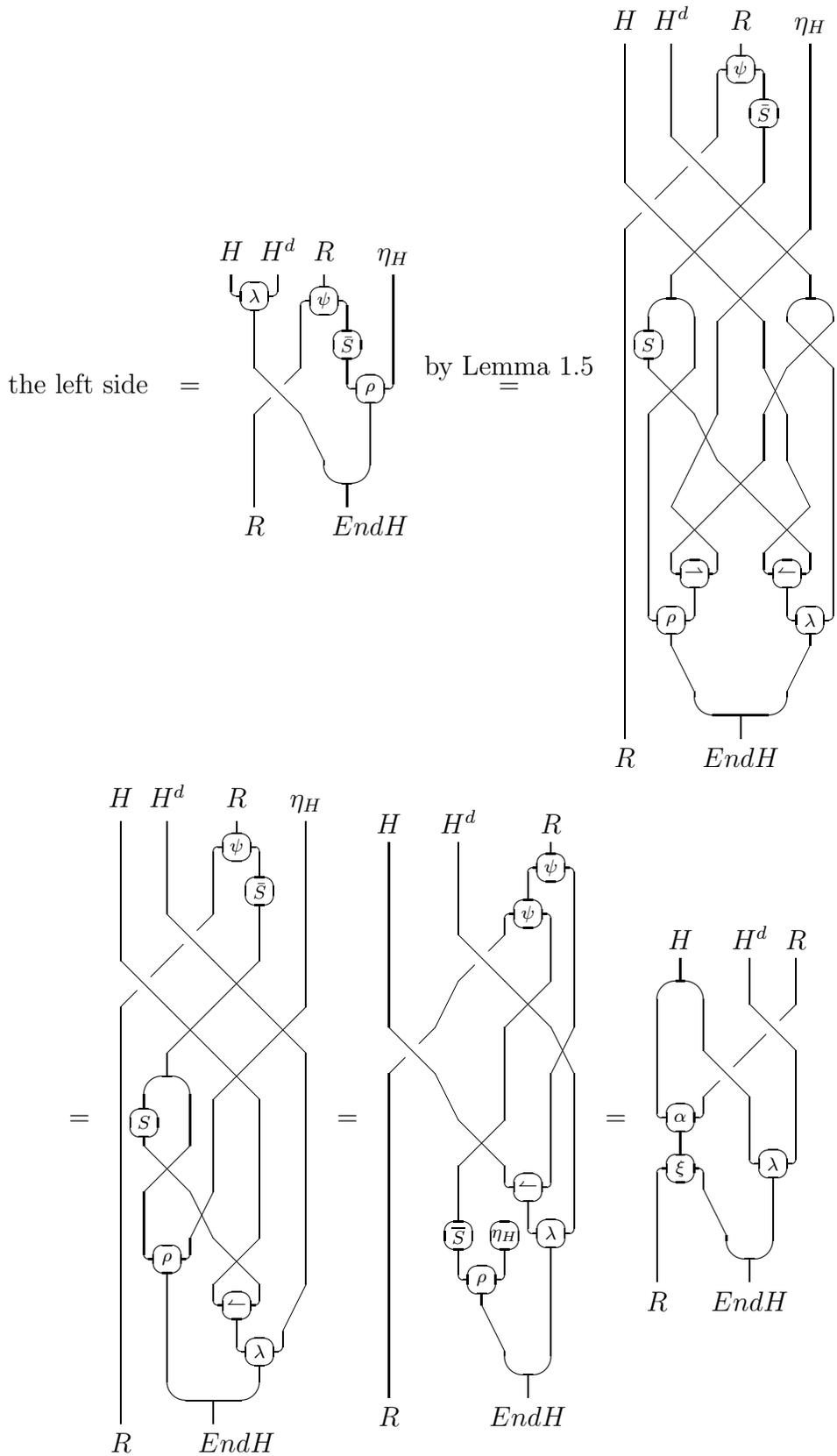
\[
\begin{tangle}
\object{H}\step[2]\object{H^{d}}\Step\object{R}\\
\tu \lambda\step\td \xi\\
\step\x\Step\id\\
\step\id\Step\cu\\
\object{R}\step[5]\object{End H}\\
\end{tangle}
\step=\step
\begin{tangle}
\step\object{H}\step[3]\object{H^{d}}\step[2]\object{R}\\
\cd\Step\x\\
\id\Step\x\Step\id\\
\tu \alpha\Step\tu \lambda\\
\td \xi\Step\dd\\
\id\Step\cu\\
\object{R}\step[4]\object{End H}\\
\end{tangle}\ \ \ .
\step[5] ......(*)
\]
\[
\hbox{the left side}\step=\step
\begin{tangle}
\object{H}\step[2]\object{H^{d}}\step[2]\object{R}\step[3]\object{\eta_H}\\
\tu \lambda\step\td \psi\step[2]\id\\
\step\id\step[2]\id\step[2]\O {\bar{S}}\step[2]\id\\

\step\x\step[2]\tu \rho\\
\step\id\Step\d\Step\id\\
\step\id\step[3]\cu\\
\step\object{R}\step[5]\object{End H}\\
\end{tangle}
\step\stackrel{\hbox{by Lemma \ref {1.6}}}{\step[3]=\step[3]}\step
\begin{tangle}
\object{H}\step[2]\object{H^{d}}\step[3]\object{R}\step[3]\object{\eta_H}\\
\id\step[2]\id\step[2] \td \psi\step[2]\id\\
\id\step[2]\id\step[2]\id\step[2]\O {\bar{S}}\step[2]\id\\
\id\Step\x\Step\id\step[2]\id\\
\x\Step\XX\step[2]\id\\
\id\Step\XX\Step\XX\\
\id\step\cd\step\XX\step\cd\\
\id\step\O S\Step\id\step\id\Step\id\step\XX\\
\id\step\XX\step\id\Step\X\Step\id\\
\id\step\id\Step\X\Step\id\step\id\Step\id\\
\id\step\id\step\dd\step\XX\step\d\step\id\\
\id\step\id\step\XX\Step\XX\step\id\\
\id\step\id\step\tu \rightharpoonup\Step
\tu\leftharpoonup\step\id\\
\id\step\tu \rho\step[4]\tu \lambda\\
\id\Step\d\step[4]\dd\\
\id\step[3]\Cu\\
\object{R}\step[5]\object{End H}\\
\end{tangle} \]
\[
\step=\step
\begin{tangle}
\object{H}\step[2]\object{H^{d}}\step[3]\object{R}\step[3]\object{\eta_H}\\
\id\step[2]\id\step[2] \td \psi\step[2]\id\\
\id\step[2]\id\step[2]\id\step[2]\O {\bar S}\step[2]\id\\
\id\Step\x\Step\id\step[2]\id\\
\x\Step\XX\step[2]\id\\
\id\Step\XX\Step\XX\\
\id\step\cd\step\XX\step[2]\id\\
\id\step\O {S}\Step\id\step\id\Step\id\Step\id\\
\id\step\XX\step\id\Step\id\Step\id\\
\id\step\id\Step\X\Step\id\Step\id\\
\id\step\tu \rho\step\XX\Step\id\\
\id\Step\id\Step\tu \leftharpoonup\step\dd\\
\id\Step\id\step[3]\tu \lambda\\
\id\Step\Cu\\
\object{R}\step[5]\object{End H}\\
\end{tangle}
\step=\step
\begin{tangle}
\object{H}\step[3]\object{H^{d}}\step[4]\object{R}\\
\id\step[3]\id\step[3]\td \psi\\
\id\step[3]\id\step[2]\td\psi\step\id\\
\id\step[3]\x\step[2]\id\step\id\\
\id\step[2]\dd\step[2]\XX\step\id\\
\x\step[3]\id\step[2]\X\\
\id\step[2]\d\step[2]\id\step[2]\id\step\id\\
\id\step[3]\XX\step[2]\id\step\id\\
\id\step[3]\id\step[2]\tu\leftharpoonup\step\id\\
\id\step[3]\O {\overline{S}}\step[2]\Q {\eta_{H}}\step\tu \lambda\\
\id\step[3]\tu \rho\step[2]\id\\
\id\step[4]\d\step[2]\id\\
\id\step[5]\cu\\
\object{R}\step[6]\object{End H}
\end{tangle}
\step=\step
\begin{tangle}
\step\object{H}\step[3]\object{H^{d}}\step[2]\object{R}\\
\cd\step[2]\id\step[2]\id\\
\id\step[2]\id\step[2]\x\\
\id\step[2]\x\step[2]\id\\
\tu \alpha\step[2]\id\step[2]\id\\
\td \xi\step[2]\tu \lambda\\
\id\step[2]\d\step[2]\id\\
\id\step[3]\cu\\
\object{R}\step[4]\object{End H}
\end{tangle}\]
Thus relation (*) holds.\\
Now, we check that $\xi$ is algebraic. We see that
\[
\begin{tangle}
\object{R}\Step\object{R}\\
\O \xi\Step\O \xi\\
\cu\\
\step\object{R\otimes End H}\\
\end{tangle}
\step=\step
\begin{tangle}
\step\object{R}\step[6]\object{R}\\
\td \psi\step[4]\td \psi\\
\id\step[2]\O {\bar{S}}\step\obox 2{\eta_H}\step\id\step[2]\O {\bar{S}}\step\obox 2{\eta_H}\\
\id\Step\tu \rho\step\dd\Step\tu \rho\\
\d\Step\x\step[4]\id\\
\step\cu\Step\Cu\\
\Step\object{R}\step[5]\object{End H}\\
\end{tangle}
\step\stackrel{\hbox{since $\rho(\overline{S}\otimes \eta_{H})$ is
algebraic}}{=}\]\[\step
\begin{tangle}
\step\object{R}\step[4]\object{R}\\
\td \psi\Step\td \psi\\
\id\Step\x\Step\id\\
\cu\step[2]\cu\\
\step\id\step[4]\O {\bar{S}}\step\obox 2{\eta_H}\\
\step\id\step[4]\tu \rho\\
\step\object{R}\step[5]\object{End H}\\
\end{tangle}
\step=\step
\begin{tangle}
\object{R}\Step\object{R}\\
\cu\\
\step\O \xi\\
\step[2]\object{R\otimes End H}\\
\end{tangle}
\step \hbox{and  obviously} \step[2]
\begin{tangle}
\object{\eta_R}\\
\id\step\\
\O \xi\\
\id\step\\
\object{R\otimes End H}\\
\end{tangle}
\step[3]=\step[3]\begin{tangle}
\object{\eta_{R\otimes End H}}\\
\id\\
\id\\
\id\\
\object{R\otimes End H}\\
\end{tangle}\ \ \ \ \ \ \ \ \ .
\]
Thus $\xi$ is algebraic.

Now we show that $\Phi'$ is algebraic.
\[
\begin{tangle}
\object{(R\# H)\# H^{d}}\step[7]\object{(R\# H)\# H^{d}}\\
\step\id\step[4]\id\\
\obox 2{\Phi'}\step[2]\obox 2{\Phi'}\\
\step\Cu\\
\step[3]\object{R\otimes End H}\\
\end{tangle}
\step=\Step
\begin{tangle}
\step\object{R}\step[2]\object{H}\step[2]\object{H^{d}}
\step[2]\object{R}\step[2]\object{H}\step[2]\object{H^{d}}\\
\td \xi\step\id\Step\id\step\td \xi\step\id\Step\id\\
\id\Step\id\step\tu \lambda\dd\Step\id\step\tu \lambda\\
\id\Step\cu\step\id\step[3]\cu\\
\d\step[2]\x\step[4]\id\\
\step\cu\step[2]\Cu\\
\step[2]\object{R}\step[6]\object{End H}\\
\end{tangle}
\step=\step
\begin{tangle}
\step\object{R}\step[2]\object{H}\step[2]\object{H^{d}}
\step[2]\object{R}\step[2]\object{H}\step[2]\object{H^{d}}\\
\td \xi\step\id\Step\id\step\td \xi\step\id\Step\id\\
\id\Step\id\step\tu \lambda\step\id\step[2]\id\step\tu \lambda\\
\id\Step\id\step[2]\x\Step\id\step[2]\id\\
\id\Step\id\Step\id\Step\cu\step\dd\\
\id\Step\x\step[3]\cu\\
\cu\Step\Cu\\
\step\object{R}\step[6]\object{End H}\\
\end{tangle}
\step\stackrel{\hbox{by (*)}}{=}\]
\[
\begin{tangle}
\step\object{R}\step[3]\object{H}\step[2]\object{H^{d}}
\step[2]\object{R}\step[2]\object{H}\step[3]\object{H^{d}}\\
\td \xi\step\cd\step\x\Step\id\Step\id\\
\id\Step\id\step\id\Step\hx\Step\id\Step\tu \lambda\\
\id\Step\id\step\tu \alpha\step\tu \lambda\step[3]\id\\
\id\Step\id\step\td \xi\Step\id\step[4]\id\\
\id\Step\hx\Step\cu\step[3]\dd\\
\cu\step\d\Step\Cu\\
\step\id\step[3]\Cu\\
\step\object{R}\step[5]\object{End H}\\
\end{tangle}
\stackrel{\hbox{by Lemma  \ref  {1.5}}}{=}
\begin{tangle}
\step\object{R}\step[3]\object{H}\step[2]\object{H^{d}}
\step[2]\object{R}\step[2]\object{H}\step[3]\object{H^{d}}\\
\td \xi\step\cd\step\x\Step\id\Step\id\\
\id\Step\id\step\id\Step\hx\step\cd\step\id\Step\id\\
\id\Step\id\step\tu \alpha\step\id\step\id\Step\X\Step\id\\
\id\Step\id\step\td \xi\step\id\step\tu \rightharpoonup\step\cu\\
\id\Step\hx\Step\id\step\cu\step[2]\ne2\\
\cu\step\id\Step\d\step\tu \lambda\\
\step\id\Step\id\step[3]\cu\\
\step\id\Step\Cu\\
\step\object{R}\step[5]\object{End H}\\
\end{tangle}
\stackrel{\hbox{ since } \xi \hbox{ is algebraic }}{=}\]
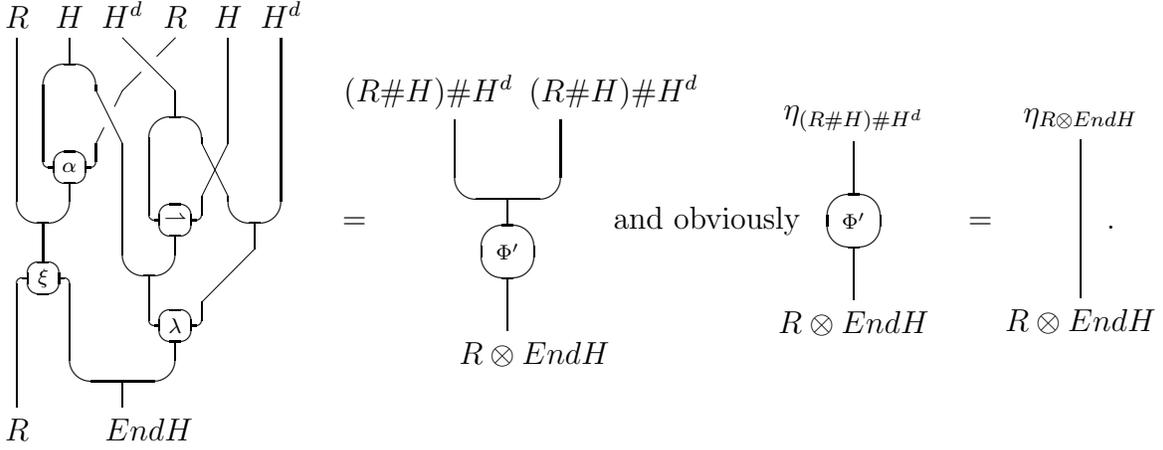
\[
\begin{tangle}
\object{R}\step[2]\object{H}\step[2]\object{H^{d}}
\step[2]\object{R}\step[2]\object{H}\step[2]\object{H^{d}}\\
\id\step\cd\step\x\Step\id\Step\id\\
\id\step\id\Step\hx\step\cd\step\id\Step\id\\
\id\step\tu \alpha\step\id\step\id\Step\X\Step\id\\
\cu\Step\id\step\tu \rightharpoonup\step\cu\\
\td \xi\Step\cu\Step\ne2\\
\id\Step\id\step[3]\tu \lambda\\
\id\Step\Cu\\
\object{R}\step[5]\object{End H}\\
\end{tangle}
\step[2]= \step[2]
\begin{tangle}
\object{(R\# H)\# H^{d}}\step[7]\object{(R\# H)\# H^{d}}\\
\step\id\step[4]\id\\
\step\Cu\\
\step[2]\obox 2{\Phi'}\\
\step[3]\id\\
\step[4]\object{R\otimes End H}\\
\end{tangle} \hbox{and obviously} \step
\begin{tangle}
\step\object{\eta_{(R\# H)\# H^{d}}}\\
\step\id\\
\obox 2{\Phi'}\\
\step\id\\
\step\object{R\otimes End H}\\
\end{tangle}
\step[3]=\step[3]
\begin{tangle}
\object{\eta_{R\otimes End H}}\\
\id\\
\id\\
\id\\
\object{R\otimes End H}\\
\end{tangle}\step .\]  Thus $\Phi '$ is algebraic. \ \
\begin{picture}(5,5)
\put(0,0){\line(0,1){5}}\put(5,5){\line(0,-1){5}}
\put(0,0){\line(1,0){5}}\put(5,5){\line(-1,0){5}}
\end{picture}\\
We obtain the following by  Theorem \ref {1.8}:
 \begin {Corollary} \label {1.9}
 Let  $H$ be a finite braided  Hopf algebra with a left dual $H^*$.
If the braiding is symmetric on $H$, then
  $$   (R \# H)\# H^{*}   \cong R \otimes (H   \# H^{*} ) \hbox { \ \ \  as }  \hbox { algebras  in } {\cal D}. $$
\end {Corollary}

This corollary reproduces   the main result in \cite {Zh03}.

\section {Hom functor in braided Yetter-Drinfeld module category} \label {e2}
In this section, we prove that if $V, W$ are in
$^{B}_{B}\mathcal{YD}(\mathcal{C})$, then $Hom(V,W)$ is also in
$^{B}_{B}\mathcal{YD}(\mathcal{C})$.

 Let $B$ be a Hopf algebra in braided tensor
category $({\mathcal C}, ^{\mathcal{C}}C)$, and $(M, \alpha, \phi )$
be a left $B$-module and left $B$ -comodule in $({\mathcal C},
^{\mathcal{C}}C)$. If

 (YD):
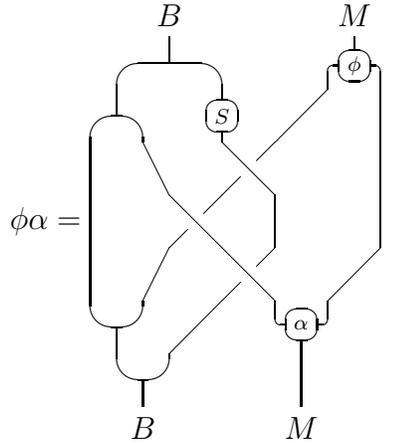
\[ \phi \alpha =  \begin{tangle}
\step [3]\object{B}\step[7]\object{M} \\

\step \Cd \step [4] \td \phi \\

\cd \step [3] \O S \step [3] \ne2 \step [2] \id \\

\id \step [2] \nw1  \step [2]  \x \step [4] \id  \\

\id \step [2]  \step \x  \step[2] \id \step [4] \id\\

\id \step [2] \ne1  \step [2] \x \step[3]   \ne2\\

\cu  \step[2] \ne2  \step[2] \tu \alpha  \\

\step \cu  \step[4]  \step \id  \\
\step [2]\object{B}\step[6]\object{M}
  \end{tangle} \ \ \ .
  \]
then $(M, \alpha, \phi )$ is called  a braided  Yetter-Drinfeld
$B$-module in ${\mathcal C}$, written as braided YD $B$-module in
short. Let $^B_B {\mathcal{ YD}(\mathcal{C})}$ denote the category
of all braided Yetter-Drinfeld $B$-modules in ${\mathcal C}$.

 Throughout this section, the braiding is in $\mathcal{C}$ and is symmetric on set $\{B,V, W, Hom(V,
 W)\}$, where $B,V,W\in \mathcal{C}$. Assume that
 $\otimes_{\mathcal{D}}=:\otimes_{\mathcal{C}}$.
 $H$ and $H^{d}$ be
braided Hopf algebra in $^B_B{\cal YD}(\mathcal{C})$ and $H^d$ is a
quasi-dual of $H$ under the operations: \\
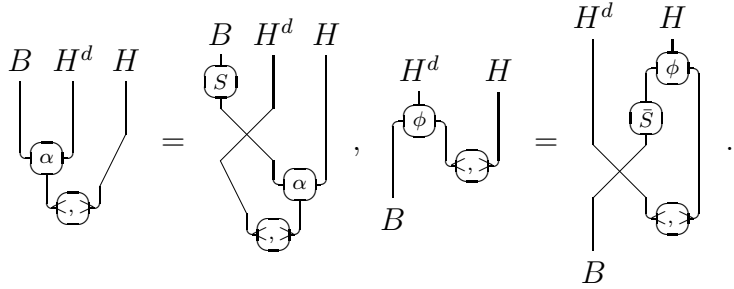
\[
\begin{tangle}
\object{B}\step[2]\object{H^{d}}\step[2]\object{H}\\
\id\step[2]\id\step[2]\id\\
\tu \alpha\step\dd\\
\step\coro {<,>}
\end{tangle}
\step = \step
\begin{tangle}
\object{B}\step[2]\object{H^{d}}\step[2]\object{H}\\
\O S\step[2]\id\step[2]\id\\
\XX\step[2]\id\\
\d\step\tu \alpha\\
\step\coro {<,>}
\end{tangle}
\step , \step
\begin{tangle}
\step\object{H^{d}}\step[3]\object{H}\\
\td \phi\step[2]\id\\
\id\step[2]\coro {<,>}\\
\object{B}
\end{tangle}
\step = \step
\begin{tangle}
\object{H^{d}}\step[3]\object{H}\\
\id\step[2]\td \phi\\
\id\step[2]\O {\bar{S}}\step[2]\id\\
\XX\step[2]\id\\
\id\step[2]\coro {<,>}\\
\object{B}
\end{tangle}\step.\]
\begin {Lemma} \label {2.1}
 If $B$ has left duals and invertible antipode, then

(i)\ \ If $(V, \alpha_V, \phi _V)$ and $(W, \alpha_W, \phi _W)$ are
in $^{B}_{B}\mathcal{YD}(\mathcal{C})$, then
 $Hom _\mathcal{C} (V, W)$ is in $^{B}_{B}\mathcal{YD}(\mathcal{C})$ under the following module operation and comodule operation:

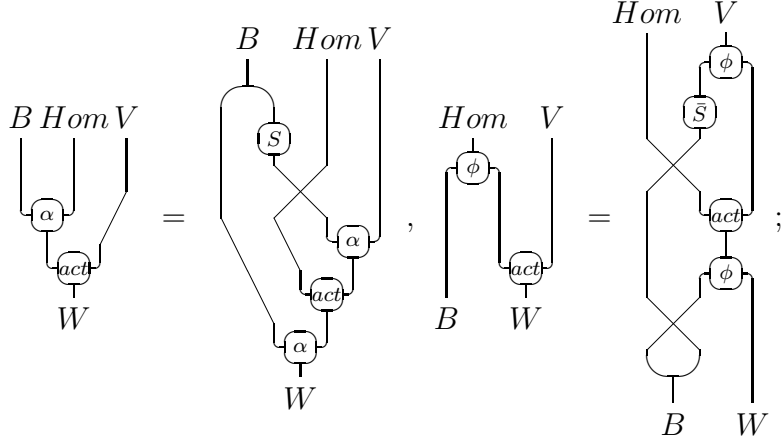
\[
\begin{tangle}
\object{B}\step[2]\object{Hom}\step[2]\object{V}\\
\id\step[2]\id\step[2]\id\\
\tu \alpha\step\dd\\
\step\tu {act}\\
\step[2]\object{W}
\end{tangle}
\step = \step
\begin{tangle}
\step\object{B}\step[3]\object{Hom}\step[2]\object{V}\\
\cd\step[2]\id\step[2]\id\\
\id\step[2]\O S\step[2]\id\step[2]\id\\
\id\step[2]\XX\step[2]\id\\
\d\step\d\step\tu \alpha\\

\step\d\step\tu {act}\\

\step[2]\tu \alpha\\
\step[3]\object{W}
\end{tangle}
\step , \step
\begin{tangle}
\step\object{Hom}\step[3]\object{V}\\
\td \phi\step[2]\id\\
\id\step[2]\id\step[2]\id\\
\id\step[2]\tu {act}\\
\object{B}\step[3]\object{W}
\end{tangle}
\step = \step
\begin{tangle}
\object{Hom}\step[3]\object{V}\\
\id\step[2]\td \phi\\
\id\step[2]\O {\bar{S}}\step[2]\id\\
\XX\step[2]\id\\
\id\step[2]\tu {act}\\
\id\step[2]\td \phi\\
\XX\step[2]\id\\
\cu\step[2]\id\\
\step\object{B}\step[3]\object{W}
\end{tangle}\ \ ;\]

(ii)\ \  $End M$ is an algebra in $^B_B{\cal YD}(\cal C),$  where
$M$ is in $^B_B{\cal YD}(\cal C)$.

\end {Lemma}
{\bf Proof.} Let $B^{*}$ denote the left deal of $B$, the
definitions of $\alpha$ and $\phi$ are reasonable, since $act$
satisfy elimination. In fact, let

\[
\begin{tangle}
\object{B^{*}}\step[2]\object{Hom}\step[2]\object{V}\\
\id\step[2]\id\step[2]\id\\
\tu {\hat{\alpha}}\step\dd\\
\step\tu {act}\\
\step[2]\object{W}
\end{tangle}
 = \step
\begin{tangle}
\object{B^{*}}\step[2.5]\object{Hom}\step[2.5]\object{V}\\
\id\step[2]\id\step[2]\td \phi\\
\id\step[2]\XX\step[2]\id\\
\id\step[2]\id \step[2]\tu {act}\\
\id\step[2]\id\step[2]\td \phi\\
\id\step[2]\id\step[2]\O S\step[2]\id\\
\d\step\cu\step[2]\id\\
\step\ev\step[3]\object{W}\\
\end{tangle}\ \ ,\step
\hbox{we can define}\step
\begin{tangle}
\step\object{Hom}\step[3]\object{V}\\
\td \phi\step[2]\id\\
\id\step[2]\id\step[2]\id\\
\id\step[2]\tu {act}\\
\object{B}\step[3]\object{W}
\end{tangle}
\step=\step
\begin{tangle}
\step[4]\object{Hom}\step[3]\object{V}\\
\coev\step[2]\id\step[2]\dd\\
\O {\bar{S}}\step[2]\tu {\hat{\alpha}}\step\dd\\
\id\step[3]\tu {act}\\
\object{B}\step[4]\object{W}
\end{tangle}
= \step
\begin{tangle}
\object{Hom}\step[3]\object{V}\\
\id\step[2]\td \phi\\
\id\step[2]\O {\bar{S}}\step[2]\id\\
\XX\step[2]\id\\
\id\step[2]\tu {act}\\
\id\step[2]\td \phi\\
\XX\step[2]\id\\
\cu\step[2]\id\\
\step\object{B}\step[3]\object{W}
\end{tangle}\]
Now we show that $Hom_{\mathcal{C}}(V,W)$ is a module and a
comodule.
\[
\begin{tangle}
\object{B}\step[2]\object{B}\step[2]\object{Hom}\step[2]\object{V}\\
\d\step\tu \alpha\step\dd\\
\step\tu \alpha\step\dd\\
\step[2]\tu {act}\\
 \step[3]\object{W}
\end{tangle}
= \step
\begin{tangle}
\step\object{B}\step[3]\object{B}\step[3]\object{Hom}\step[2]\object{V}\\
\cd\step\cd\step[2]\id\step[2]\id\\
\id\step[2]\O S\step\id\step[2]\O S\step[2]\id\step[2]\id\\
\id\step[2]\X\step[2]\id\step[2]\id\step[2]\id\\
\id\step[2]\id\step\XX\step[2]\id\step[2]\id\\
\id\step[2]\id\step\id\step[2]\XX\step[2]\id\\
\id\step[2]\id\step\id\step[2]\id\step[2]\tu \alpha\\
\id\step[2]\id\step\XX\step[2]\dd\\
\id\step[2]\id\step\d\step\tu \alpha\\
\d\step\d\step\tu {act}\\
\step\d\step\tu \alpha\\
\step[2]\tu \alpha\\
 \step[3]\object{W}
\end{tangle}
\step =
\begin{tangle}
\step\object{B}\step[3]\object{B}\step[3]\object{Hom}\step[2]\object{V}\\
\cd\step\cd\step[2]\id\step[2]\id\\
\id\step[2]\O S\step\id\step[2]\O S\step[2]\id\step[2]\id\\
\id\step[2]\X\step[2]\id\step[2]\id\step[2]\id\\
\cu\step\XX\step[2]\id\step[2]\id\\
\step\id\step[2]\id\step[2]\XX\step[2]\id\\
\step\id\step[2]\XX\step[2]\id\step[2]\id\\
\step\d\step\d\step\cu\step\dd\\
\step[2]\d\step\d\step\tu \alpha\\
\step[3]\d\step\tu {act}\\
\step[4]\tu \alpha\\
 \step[5]\object{W}
\end{tangle}
\step =\step
\begin{tangle}
\step\object{B}\step[3]\object{B}\step[3]\object{Hom}\step[2]\object{V}\\
\cd\step\cd\step[2]\id\step[2]\id\\
\id\step[2]\X\step[2]\id\step[2]\id\step[2]\id\\
\cu\step\cu\step[2]\id\step[2]\id\\
\step\id\step[3]\O S\step[2]\dd\step[2]\id\\
\step\d\step[2]\XX\step[2]\dd\\
\step[2]\d\step\d\step\tu \alpha\\
\step[3]\d\step\tu {act}\\
\step[4]\tu \alpha\\
 \step[5]\object{W}
\end{tangle}
\step\]\[= \step
 \begin{tangle}
\object{B}\step[2]\object{B}\step[2]\object{Hom}\step[2]\object{V}\\
\cu\step[2]\id\step[2]\id\\
\cd\step[2]\id\step[2]\id\\
\id\step[2]\O S\step[2]\id\step[2]\id\\
\id\step[2]\XX\step[2]\id\\
\d\step\d\step\tu \alpha\\
\step\d\step\tu {act}\\
\step[2]\tu \alpha\\
\step[3]\object{W}
\end{tangle}\step =
\begin{tangle}
\object{B}\step[2]\object{B}\step[2]\object{Hom}\step[2]\object{V}\\
\cu\step\dd\step\dd\\
\step\tu \alpha\step\dd\\
\step[2]\tu {act}\\
\step[3]\object{W}
\end{tangle};\step
\begin{tangle}
\step\object{Hom}\step[4]\object{V}\\
\td \phi\step[3]\id\\
\id\step\td \phi\step[2]\id\\
\id\step\id\step[2]\tu {act}\\
\object{B}\step\object{B}\step[3]\object{W}
\end{tangle}
\step\]\[ = \step
\begin{tangle}
\object{Hom}\step[4]\object{V}\\
\id \step[3]\td \phi\\
\d\step[2]\O {\bar{S}}\step\td \phi\\
\step\XX\step\O {\bar{S}}\step[2]\id\\
\dd\step[2]\X\step[2]\id\\
\id\step[2]\dd \step\tu {act}\\
\id\step[2]\id \step[2]\td \phi\\
\id\step[2]\XX\step\td \phi\\
\d\step\cu\step\id\step[2]\id\\
\step\XX\step[2]\id\step[2]\id\\
\step\id\step[2]\XX\step[2]\id\\
\step\id\step[2]\cu\step[2]\id\\
\step\object{B}\step[3]\object{B}\step[3]\object{W}
\end{tangle}
\stackrel{by \cite [Proposition 1.0.13 ] {Zh99} }{=}
\begin{tangle}
\step\object{Hom}\step[4]\object{V}\\
\step\id\step[3]\td \phi\\
\step\d\step[2]\O {\bar{S}}\step[2]\id\\
\step[2]\XX\step[2]\id\\
\step\dd\step[2]\tu {act}\\
\step\id\step[3]\td \phi\\
\cd\step\cd\step\id\\
\XX\step\id\step[2]\id\step\id\\
\id\step[2]\X\step[2]\id\step\id\\
\id\step[2]\hddcu\step[2]\id\step\id\\
\XX\step[2]\dd\step\id\\
\id\step[2]\XX\step[2]\id\\
\id\step[2]\cu\step[2]\id\\
 \object{B}\step[3]\object{B}\step[3]\object{W}
\end{tangle}
\step =
\begin{tangle}
\step\object{Hom}\step[4]\object{V}\\
\step\id\step[3]\td \phi\\
\step\d\step[2]\O {\bar{S}}\step[2]\id\\
\step[2]\XX\step[2]\id\\
\step\dd\step[2]\tu {act}\\
\step\d\step[2]\td \phi\\
\step[2]\XX\step[2]\id\\
\step\dd\step[2]\id\step[2]\id\\
\cd\step\cd\step\id\\
\id\step[2]\X\step[2]\id\step\id\\
\cu\step\cu\step\id\\
 \step\object{B}\step[3]\object{B}\step[2]\object{W}
\end{tangle}
\step =
\begin{tangle}
\step[2]\object{Hom}\step[3]\object{V}\\
\step\td \phi\step[2]\id\\
\step\id\step[2]\id\step[2]\id\\
\step\id\step[2]\tu {act}\\
\cd\step[2]\id\\
 \object{B}\step[2]\object{B}\step[2]\object{W}
\end{tangle}\ \ .
\]
Now we show that
\[ \phi \alpha \step= \step\begin{tangle}
\step [3]\object{B}\step[5]\object{Hom} \\
\step \Cd \step [2] \td \phi \\
\cd \step [3] \O S \step [2] \id \step [2] \id \\
\id \step [2] \nw1  \step [2]  \XX \step [2] \id  \\
\id \step [2]  \step \XX  \step[2] \id \step [2] \id\\
\id \step [2] \ne1  \step [2] \XX \step[2]   \id\\
\cu  \step[2] \ne2  \step[2] \tu \alpha  \\
\step \cu  \step[4]  \step \id  \\
\step [2]\object{B}\step[6]\object{Hom}
\end{tangle}
\step=\step
\begin{tangle}
\step [3]\object{B}\step[6]\object{Hom} \\
\step \Cd \step [3] \td \phi \\
\step\id\step[3]\cd\step[2]\id\step[2]\id\\
\step\id\step[3]\id\step[2]\O S\step[2]\id\step[2]\id\\
\step\id\step[3]\id\step[2]\XX\step[2]\id\\
\step\id\step[3]\d\step\cu\step[2]\id\\
\step\id\step[4]\XX\step[2]\dd\\
\step\Cu\step[2]\tu \alpha\\
\step [3]\object{B}\step[5]\object{Hom}
\end{tangle}\step.\] In  fact ,
\[\begin{tangle}
\step [3]\object{B}\step[6]\object{Hom}\step[2]\object{V} \\
\step \Cd \step [3] \td \phi \step\id\\
\step\id\step[3]\cd\step[2]\id\step[2]\id\step\id\\
\step\id\step[3]\id\step[2]\O S\step[2]\id\step[2]\id\step\id\\
\step\id\step[3]\id\step[2]\XX\step[2]\id\step\id\\
\step\id\step[3]\d\step\cu\step[2]\id\step\id\\
\step\id\step[4]\XX\step[2]\dd\step\id\\
\step\Cu\step[2]\tu \alpha\step\dd\\
\step[3]\id\step[5]\tu {act}\\
 \step [3]\object{B}\step[5]\step\object{W}
\end{tangle}
\step=\]\[
\begin{tangle}
\step [2]\object{B}\step[6]\object{Hom}\step[3]\object{V} \\
\Cd\step[4]\id\step[2]\td \phi\\
\id\step[3]\cd\step[3]\XX\step[2]\id\\
\id\step[2]\cd\step\O S\step[3]\id\step[2]\id\step[2]\id\\
\id\step[2]\id\step[2]\O
S\step\id\step[3]\id\step[2]\id\step[2]\id\\
\id\step[2]\id\step[2]\X
\step[3]\id\step[2]\id\step[2]\id\\
\id\step[2]\id\step\dd\cd\step[2]\id\step[2]\id\step[2]\id\\
\id\step[2]\id\step\id\step\id\step\cd\step\id\step[2]\id\step[2]\id\\
\id\step[2]\id\step\id\step\id\step\id\step[2]\O S\step\id\step[2]\id\step[2]\id\\
\id\step[2]\id\step\id\step\id\step\id\step[2]\X\step[2]\id\step[2]\id\\
\id\step[2]\id\step\id\step\id\step\d\step\hdcu\step[2]\id\step[2]\id\\
\id\step[2]\id\step\id\step\id\step[2]\XX\step[2]\id\step[2]\id\\
\id\step[2]\id\step\id\step\cu\step[2]\id\step[2]\id\step[2]\id\\
\id\step[2]\id\step\id\step[2]\O {\bar{S}}\step[3]\XX\step[2]\id\\
\id\step[2]\id\step\id\step[2]\d\step[2]\d\step\tu
\alpha\\
\id\step[2]\id\step\d\step[2]\d\step[2]\tu
{act}\\
\id\step[2]\d\step\d\step[2]\id\step[2]\td
\phi\\
\d\step[2]\d\step\id\step[2]\XX\step[2]\id\\
\step\d\step[2]\id\step\d\step\cu\step[2]\id\\
\step[2]\d\step\d\step\cu\step[2]\dd\\
\step[3]\id\step[2]\XX\step[2]\dd\\
\step[3]\cu\step[2]\tu \alpha\\
 \step [4]\object{B}\step[4]\object{W}
\end{tangle}
\stackrel{by \cite [Proposition 1.0.13 ]{Zh99}}{=}
\begin{tangle}
\step \object{B}\step[4]\object{Hom}\step[3]\object{V} \\
\cd\step[3]\id\step[2]\td \phi\\
\id\step\cd\step[2]\id\step[2]\O {\bar{S}}\step[2]\id\\
\id\step\id\step\cd\step\XX\step[2]\id\\
\id\step\id\step\XX\step\id\step[2]\id\step[2]\id\\
\id\step\id\step\O S\step[2]\O S\step\id\step[2]\id\step[2]\id\\
\id\step\id\step\id\step[2]\hddcu\step[2]\id\step[2]\id\\
\id\step\id\step\XX\step[2]\dd\step[2]\id\\
\id\step\id\step\id\step[2]\XX\step[2]\dd\\
\id\step\id\step\id\step[2]\d\step\tu \alpha\\
\id\step\id\step\id\step[3]\tu {act}\\
\id\step\id\step\d\step[2]\td \phi\\
\id\step\id\step[2]\XX\step[2]\id\\
\id\step\d\step\cu\step[2]\id\\
\id\step[2]\XX\step[2]\dd\\
\cu\step[2]\tu \alpha\\
 \step \object{B}\step[4]\object{W}
\end{tangle}
= \begin{tangle}
\step[2] \object{B}\step[5]\object{Hom}\step[3]\object{V} \\
\Cd\step[3]\id\step[2]\td \phi\\
\id\step[3]\cd\step[2]\id\step[2]\O {\bar{S}}\step[2]\id\\
\id\step[2]\cd\step\O S\step[2]\XX\step[2]\id\\
\id\step[2]\id\step[2]\O
S\step\id\step[2]\id\step[2]\id\step[2]\id\\
\id\step[2]\id\step[2]\id\step\XX\step[2]\id\step[2]\id\\
\id\step[2]\id\step[2]\id\step\id\step[2]\XX\step[2]\id\\
\id\step[2]\id\step[2]\id\step\d\step\d\step\tu \alpha\\
\id\step[2]\id\step[2]\cu\step[2]\tu {act}\\
\d\step\d\step[2]\id\step[3]\td \phi\\
\step\d\step\d\step\d\step[2]\id\step[2]\id\\
\step[2]\d\step\d\step\XX\step[2]\id\\
\step[3]\id\step[2]\id\step\cu\step[2]\id\\
\step[3]\id\step[2]\XX\step[2]\dd\\
\step[3]\cu\step[2]\tu \alpha\\
 \step [4]\object{B}\step[4]\object{W}
\end{tangle},\]\[ \step
\begin{tangle}
\object{B}\step[2]\object{Hom}\step[2]\object{V} \\
\tu \alpha\step [2] \id\\
\td \phi\step [2] \id \\
\id \step [2] \tu {act} \\
\object{B}\step[3]\object{W}
\end{tangle}\step=\]
\[
\begin{tangle}
\step\object{B}\step[3]\object{Hom}\step[3]\object{V} \\
\cd\step[2]\id\step[2]\td \phi\\
\id\step[2]\O S\step[2]\id\step[2]\O {\bar{S}}\step[2]\id\\
\id\step[2]\id\step[2]\XX\step[2]\id\\
\id\step[2]\XX\step[2]\id\step[2]\id\\
\d\step\d\step\XX\step[2]\id\\
\step\d\step\d\d\step\tu \alpha\\
\step[2]\XX\step\tu {act}\\
\step[2]\id\step[2]\tu \alpha\\
\step[2]\id\step[2]\td \phi\\
\step[2]\XX\step[2]\id\\
\step[2]\cu\step[2]\id\\
\step[3]\object{B}\step[3]\object{W}
\end{tangle}\step=\step
\begin{tangle}
\step [2]\object{B}\step[4]\object{Hom}\step[3]\object{V} \\
\Cd\step[2]\id\step[2]\td \phi\\
\id\step[4]\O S\step[2]\id\step[2]\O {\bar{S}}\step[2]\id\\
\d\step[3]\id\step[2]\XX\step[2]\id\\
\step\d\step[2]\XX\step[2]\id\step[2]\id\\
\step[2]\XX\step[2]\XX\step[2]\id\\
\step[2]\id\step\cd\step\d\step\tu \alpha\\
\step[2]\id\step\id\step\cd\step\tu {act}\\
\step[2]\id\step\id\step\id\step[2]\id\step\td \phi\\
\step[2]\id\step\id\step\id\step[2]\O S\step\id\step[2]\id\\
\step[2]\id\step\id\step\id\step[2]\X\step[2]\id\\
\step[2]\id\step\id\step\d\step\hddcu\step[2]\id\\
\step[2]\id\step\id\step[2]\X\step[2]\dd\\
\step[2]\id\step\cu\step\tu \alpha\\
\step[2]\XX\step[3]\id\\
\step[2]\cu\step[3]\id\\
\step[3] \object{B}\step[4]\object{W}
\end{tangle}\step
\step=\step
\begin{tangle}
\step [2]\object{B}\step[5]\object{Hom}\step[4]\object{V} \\
\Cd\step[3]\id\step[2]\td \phi\\
\id\step[3]\cd\step[2]\id\step[2]\O {\bar{S}}\step[2]\id\\
\id\step[2]\cd\step\O S\step[2]\XX\step[2]\id\\
\id\step[2]\id\step[2]\O S\step\XX\step[2]\id\step[2]\id\\
\id\step[2]\id\step[2]\id\step\id\step[2]\XX\step[2]\id\\
\id\step[2]\id\step[2]\id\step\id\step[2]\d\step\tu \alpha\\
\id\step[2]\id\step[2]\hdcu\step[3]\tu {act}\\
\d\step\d\step[2]\d\step[2]\td \phi\\
\step\d\step\d\step[2]\XX\step[2]\id\\
\step[2]\d\step\d\step\cu\step[2]\id\\
\step[3]\d\step\XX\step[2]\dd\\
\step[4]\hcu\step[2]\tu \alpha\\
 \step [5]\object{B}\step[4]\object{W}
\end{tangle}
\]
So $Hom _\mathcal{C} (V, W)$ is in
$^{B}_{B}\mathcal{YD}(\mathcal{C})$.

(ii) It is straightforward.\ \ $\Box$

\begin{Lemma}
Let $H$ is a braided Hopf algebra in
$^{B}_{B}\mathcal{YD}(\mathcal{C})$, $E=: End H$, then

(i)$act: E\otimes H\rightarrow H$ is a morphism in
$^{B}_{B}\mathcal{YD}(\mathcal{C})$.

(ii)The evaluation $<,>$ is in $^{B}_{B}\mathcal{YD}(\mathcal{C})$.

\end{Lemma}
{\bf Proof}
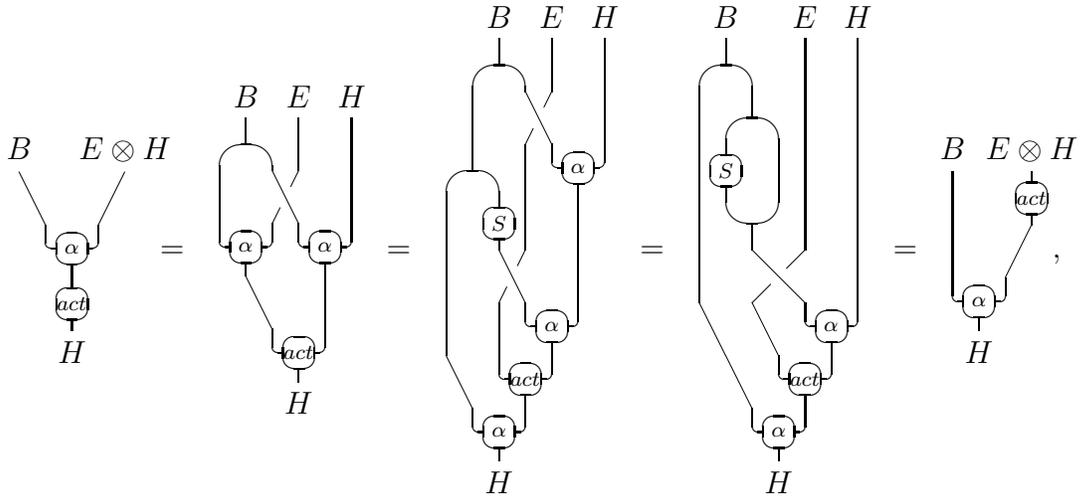
\[\begin{tangle}
\object{B}\step[4]\object{E\otimes H}\\
\d\step[2]\dd\\
\step\tu \alpha\\
\step[2]\O {act}\\
\step[2]\object{H}
\end{tangle}
\step=\step
\begin{tangle}
\step\object{B}\step[2]\object{E}\step[2]\object{H}\\
\cd\step\id\step[2]\id\\
\id\step[2]\hx\step[2]\id\\
\tu \alpha\step\tu \alpha\\
\step\d\step[2]\id\\
\step[2]\tu {act}\\
\step[3]\object{H}
\end{tangle}
\step=\step
\begin{tangle}
\step[2]\object{B}\step[2]\object{E}\step[2]\object{H}\\
\step\cd\step\id\step[2]\id\\
\step\id\step[2]\hx\step[2]\id\\
\cd\step\id\step\tu \alpha\\
\id\step[2]\O S\step\id\step[2]\id\\
\id\step[2]\hx\step[2]\id\\
\id\step[2]\id\step\tu \alpha\\
\d\step\tu {act}\\
\step\tu \alpha\\
\step[2]\object{H}
\end{tangle}
\step=\step
\begin{tangle}
\step\object{B}\step[3]\object{E}\step[2]\object{H}\\
\cd\step[2]\id\step[2]\id\\
\id\step\cd\step\id\step[2]\id\\
\id\step\O S\step[2]\id\step\id\step[2]\id\\
\id\step\cu\step\id\step[2]\id\\
\id\step[2]\x\step[2]\id\\
\d\step\d\step\tu \alpha\\
\step\d\step\tu {act}\\
\step[2]\tu \alpha\\
 \step[3]\object{H}
\end{tangle}
\step=\step
\begin{tangle}
\object{B}\step[3]\object{E\otimes H}\\
\id\step[3]\O {act}\\
\id\step[2]\dd\\
\tu \alpha\\
\step\object{H}
\end{tangle} \ \  ,
\]
so $act$ is a $B$-module homomorphism, similarly, we can show that
$act$ is a $B$-comodule homomorphism.

(ii) It is similar to (i). $\Box$\\

 If $act$ satisfy elimination, let $\mathcal{D}=^{B}_{B}\mathcal{YD}(\mathcal{C})$, $CRL1$ could be instead of:\\
$\mathbf{CRL1^{'}}$ \ \ \ $E =: End _{\cal C} \ H$ satisfy:\\
\[
\begin{tangle}
\object{E}\step[2] \object{E} \step\object{H} \\
\tu {m} \step\id\\
\step\tu {act}\\
\step[2]\object{H}\\
\end{tangle}
\step=\step
\begin{tangle}
\object{E}\step \object{E} \step[2]\object{H}  \\
\id\step\tu {act}\\
\tu {act}\\
\step\object{H}\\
\end{tangle} \step,\step
\begin{tangle}
\step[2]\object{H}\\
\Q {\eta_E}\step[2]\id \\
\tu {act}\\
\step\object{H}\\
\end{tangle}
\step=\step
\begin{tangle}
\object{H}\\
\id \\
\id \\
\object{H}
\end{tangle}\step.
\]

 \section{Duality theorems in  Yetter-Drinfeld module categories}\label {e3}

 In this section, we present the duality theorem for braided Hopf algebras
in the Yetter-Drinfeld module category $^B_B {\cal YD}$(i.e. if
$\mathcal{C}$ is the category of vector spaces, we write $^B_B
{\mathcal {YD} (\mathcal{C})} = ^B_B {\mathcal {YD}}$). Throughout
this section, $H$ is a braided Hopf algebra in $^B_B{\cal YD}$ with
Hopf algebra $B$ and $H^d$ is a quasi-dual of $H$ under a left
faithful $<, >$ (i.e. $<x, H > = 0$ implies $x =0$)such that $(b
\cdot f)(x) = \sum b_2\cdot f(S(b_1)\cdot x )$ and $\sum <f _{(0)},
x>f_{(-1)} = \sum <f, x_{(0)}> S^{-1}(x_{(-1)}) $ for any $x\in H, b
\in B, f \in H^d.$ Let $<,
>_{ev}$ the ordinary evaluation of any spaces.

Let $A$ be any braided algebra in $\cal{D}=$  $^B_B {\cal YD}$.
Define
$$A^\circ _{\cal D}= \{f \in A^* \mid Ker (f) \hbox
{ contains an ideal of finite codimension in } {}^B_B{\cal YD} \}
$$
Consequently, let $H$ be a Hopf algebra in $^B_B {\cal YD}$, $U$ be
a subHopfalgebra of $H^\circ _{\cal D}$. Then we call that $U$
satisfy the $RL$-condition with respect to $H$ if
$\rho(U\#1)\subseteq \lambda(H\#U)$ \cite [Definition 9.4.5] {Mo93}.
\begin{Lemma}\label {3.1}
If braided algebra $A\in ob\mathcal{D}$, then
$A_{\mathcal{D}}^{\circ}\in ob\mathcal{D}$.
\end{Lemma}
{\bf Proof.}  By Lemma \ref {2.1}, $A^*\in ob(\mathcal{D})$. For any
$f \in A^\circ_{\mathcal{D}}$, there exists an ideal $I$ of $A$ and
$I$ is a $B$-submodule and a $B$-subcomodule of $A$ with finite
codimension and $f(I)=0$. Since $(b\cdot f) (x) = \sum b_{2}\cdot f
(S(b_{1})\cdot x) = 0$ for any $b \in B, x \in I$, we  have  $b\cdot
f \in A^\circ_{\mathcal{D}} $. Thus $A^\circ_{\mathcal{D}}$ is a
$B$-submodule of $A^*$. By Lemma \ref {2.1}, we can assume $\phi
_{A^*} (f) = \sum_i u_i \otimes v_i$ with linear independent
$u_i's.$ Since $\sum_i u_i v_i (x) = \sum f (x_{(0)})S^{-1}
(x_{(-1)})=0 $  for any $x\in I$, we have that $v_i (x)=0$ and
$v_i(I)=0$, which implies $v_i \in A^\circ_{\mathcal{D}}$. thus
$A^{\circ}_{\mathcal{D}}$ is a $B$-subcomodule of $A^*.$
Consequently, it is clear that $A_{\mathcal{D}}^{\circ}\in
ob\mathcal{D}$.
\begin{Lemma}\label{3.2}
If $f$ is a morphism from $U$ to $V$ in  $\mathcal{D}$, then $f^{*}$
is a morphism from $V^{*}$ to $U^{*}$ in  $\mathcal{D}$.
\end{Lemma}

{\bf Proof.}   For any $v^* \in V^*, u \in U, b\in B,$ see that
\begin {eqnarray*}
(b \cdot f ^* (v^*))(u) &=& \sum b_{2}\cdot f^* (v^*)(S(b_{1}) \cdot u)\\
&=& \sum b_{2}\cdot v^* (f (S(b_{1}) \cdot u)) \hbox { \ \ since } f
\hbox { is a } B
\hbox {-module homomorphism }\\
&=&\sum b_{2}\cdot v^{*}( S(b_{1}) \cdot f(u))\\
&=&f^{*}(b\cdot v^{*})(u).  \end {eqnarray*} Thus $b \cdot f^*
(v^*)= f^* (b\cdot v^*)$ and $f^*$ is a $B$- module homomorphism.
Similarly, we can show that $f^*$ is a $B$-comodule
homomorphism.$\Box$

By Lemma \ref {3.1}, Lemma \ref {3.2} and \cite [Theorem
9.1.3]{Mo93}, we get the following:
\begin{Theorem}\label {3.3}
If $H$ is a  Hopf algebra in $\mathcal{D}$, then
$H^{\circ}_{\mathcal{D}}$ is a Hopf algebra in $\mathcal{D}$.
\end{Theorem}
Next, we give an example which showed that there exists  a Hopf
algebra $H$ in Yetter-Drinfeld module category such that
$H^{\circ}_{\mathcal{D}}$ is nontrival.

Let $T=\mathcal{T}$$(G,g_{i},\chi_{i};J)$ be the free algebra
generated by set $X=\{x_{i}\mid i\in J\}$ where $G$ is a group,
$J=\{1,2,\cdots,\theta\}$, $g_{i}\in Z(G)$ and $\chi_{i}\in \hat{G}$
with $i\in J$. We present the construction of $T$ as follows: Denote
by $T_{0}=\emptyset$, $T_{1}=X$, and for $n\geq 2$ by
$T_{n}=X\otimes X\otimes \cdots \otimes X$, the tensor product of
$n$ copies of the set $X$, then $T=\oplus_{n\geq 0}T_{n}$.
 Define coalgebra operations and $kG$-(co-)module operations in
 $T$ as follows:
$$\Delta x_i = x_i \otimes 1 + 1 \otimes x_i, \ \ \epsilon (x_i)
=0,$$
 $$\delta ^-(x_{i}) = g_{i}\otimes x_{i}, \qquad h \cdot x_{i} =
\chi_{i}(h)x_{i}. $$ Then  $T$ is  called a quantum tensor algebra
in $\mathcal{D}=^{kG}_{kG} {\cal YD} $.

If $\chi_{i}(g_{j})\chi_{j}(g_{i})=1$ for $i,j \in \{1,2, \cdots
\theta\}$, it is easily to check that $T$ is quantum cocommutative.

For quantum tensor algebra $T$, we can find an ideal $D$ of $T^{*}$
such that $D\subseteq T^{\circ}_{\mathcal{D}}$. We know that
$T^{*}\cong \prod_{n\geq 0}T_{n}^{*}$, denote $D=\oplus_{n\geq
0}T^{*}_{n}\subseteq T^{*}$. For any $f=\sum_{i=1}^{s}f_{i}\in D$,
denote $L_{n}=\oplus_{l>n}T_{l}$ with $n\geq s$, it is easily to
prove that $L_{n}$ is a finite codimension ideal of $T$ and
$f(L_{n})=0$. Consequently, we only need to show that $L_{n}$ is a
sub-(co-)module, it is sufficient to prove that $T_{m}$ is a
sub-(co-)module for $m\geq 0$. In fact, for any $g\in G$, $y=y_{1}
y_{2} \cdots y_{m}\in T_{m}$(i.e.the multiplication of $T$ is
$m(x\otimes y)=xy$.), where $y_{i}\in X$ for any $i\in J$, then
\begin {eqnarray*}
g\cdot y&=&g\cdot(y_{1} y_{2} \cdots  y_{m})\hbox { \ \ \ \ \ \ \ \
 since } T \hbox { is a }
\hbox {kG-module algebra }\\
&=&(g\cdot y_{1})(g\cdot y_{2})(g\cdot
y_{m})\\
&=&\chi_{1}(g)\chi_{2}(g)\cdots \chi_{m}(g)y_{1} y_{2} \cdots
y_{m}\hbox { \ \ \ \ \ Let }  \alpha \hbox { = }
\chi_{1}(g)\chi_{2}(g)\cdots \chi_{m}(g)
\hbox {\ \ }\\
&=&\alpha y\in T_{m} \end {eqnarray*} and
\begin {eqnarray*}
\delta^{-}(y)&=&\delta(y_{1} y_{2} \cdots y_{m})\hbox { \ \ \ \ \ \
\ since } T \hbox { is a }
\hbox {kG-comodule algebra }\\
&=&(g_{1}\cdot y_{1})(g_{2}\cdot y_{2})\cdots(g_{m}\cdot
y_{m})\\
&=&\chi_{1}(g_{1})\chi_{2}(g_{2})\cdots \chi_{m}(g_{m})y_{1} y_{2}
\cdots y_{m}\hbox { \ \ \ \ \ Let }  \beta \hbox { =  }
\chi_{1}(g_{1})\chi_{2}(g_{2})\cdots \chi_{m}(g_{m})\hbox {\ \
}\\
&=&\beta y\in T_{m}\subseteq KG\otimes T_{m},\end {eqnarray*} so
$T_{m}$ is a sub-(co-)module for $m\geq 0$.
\begin{Definition}\label {3.8}
Let $H$ be a braided Hopf algebra in $\mathcal{D}=^B_B{\cal YD}$,
$U$ is a Hopf subalgebra of $H_{\mathcal{D}}^{\circ}$, we define
morphisms $\lambda: H\#U\rightarrow End H$ via
$\lambda(h\#f)(k)=\Sigma hk_{1}<f,k_{2}>$, and $\rho:
U\#H\rightarrow End H$ via $\rho(f\#h)(k)=\Sigma k_{2}h<f,k_{1}>$
for all $h,k\in H, f\in U$.
\end{Definition}

It is clear that $\lambda,\rho \in \mathcal{D}$, then, since $\bar
\lambda \lambda  = id _{H\#
 H^{d}}\in \mathcal{D}$(by Lemma \ref {3.5}), $\bar{\lambda}\in
 \mathcal{D}$.

\begin {Lemma} \label {3.5}
 (i) Define $act: End H\otimes H\rightarrow H$ via $act(f\otimes
h)=f(h)$, then $act$ satisfy elimination.

(ii) If the antipode of $H$ is  invertible, then there exists
morphism $\bar \lambda_{1} $
 from $Im \lambda_{1} $ to $H\# H^{d}$ such that $\bar \lambda_{1} \lambda  = id _{H\#
 H^{d}}$. Furthermore, let $E= E_{1} \oplus Im\lambda_{1}$ where $E_{1}$ is a subspace of $E$ , define $\bar{\lambda}=0+\bar{\lambda_{2}}: E\rightarrow H\#
 H^{d}$, then $\bar \lambda \lambda  = id _{H\#
 H^{d}}$.

(iii)  If $H$ is quantum  cocommutative, then $RL$-condition holds
on $H$ and $U$ under $<,>$.
\end {Lemma}
{\bf Proof.} (i) It is straightforward.

(ii) We  define a morphism $\lambda '$ and $\bar{\lambda_{2}}$ as
follows:
\[
\begin{tangle}
\object{H\#H^{d}}\step[2.5]\object{H}\\
\O {\lambda'} \step[2]\id\\
\tu {act} \\
\step\object{H}
\end{tangle}
\step=\step
\begin{tangle}
\object{H}\step[2]\object{H^{d}}\step[2]\object{H}\\
\id\step[2]\coro {<,>}\\
\object{H}
\end{tangle}\step,\step
\begin{tangle}
\object{EndH}\step[2.5]\object{H}\\
\O {\bar{\lambda_{2}}} \step[2]\id\\
\tu {act} \\
\step\object{H}
\end{tangle}
\step=\step
\begin{tangle}
\object{EndH}\step[3]\object{H}\\
\id\step[2]\cd\\
\id\step[2]\O {\bar{S}}\step[2]\id\\
\id\step[2]\XX\\
\tu {act}\step\dd\\
\step\cu\\
\step[2]\object{H}
\end{tangle}\step,\]
obviously, $\lambda'$ is a injective. Now we show that
$\bar{\lambda_{2}}\lambda=\lambda'$.
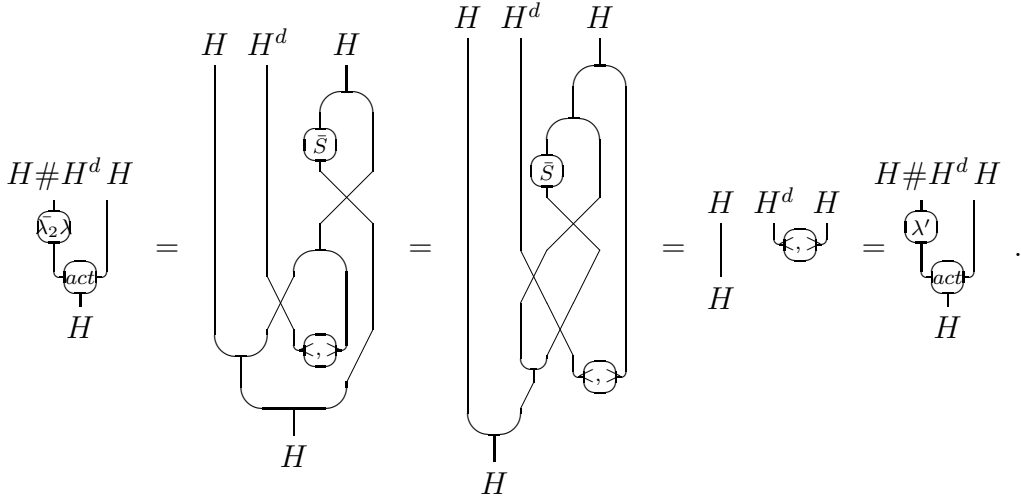
\[
\begin{tangle}
\object{H\#H^{d}}\step[2.5]\object{H}\\
\O {\bar{\lambda_{2}}\lambda} \step[2]\id\\
\tu {act} \\
\step\object{H}
\end{tangle}
\step=\step
\begin{tangle}
\object{H}\step[2]\object{H^{d}}\step[3]\object{H}\\
\id\step[2]\id\step[2]\cd\\
\id\step[2]\id\step[2]\O {\bar{S}}\step[2]\id\\
\id\step[2]\id\step[2]\XX\\
\id\step[2]\id\step[1]\cd\step\id\\
\id\step[2]\X\step[2]\id\step\id\\
\cu\step\coro {<,>}\dd\\
\step\Cu\\
\step[3]\object{H}
\end{tangle}\step=\step
\begin{tangle}
\object{H}\step[2]\object{H^{d}}\step[3]\object{H}\\
\id\step[2]\id\step[2]\cd\\
\id\step[2]\id\step\cd\step\id\\
\id\step[2]\id\step\O {\bar{S}}\step[2]\id\step\id\\
\id\step[2]\id\step\XX\step\id\\
\id\step[2]\X\step\dd\step\id\\
\id\step[2]\id\step\X\step[2]\id\\
\id\step[2]\hddcu\step\coro{<,>}\\
\cu\\
\step\object{H}
\end{tangle}\step=\step
\begin{tangle}
\object{H}\step[2]\object{H^{d}}\step[2]\object{H}\\
\id\step[2]\coro {<,>}\\
\object{H}
\end{tangle}\step=\step
\begin{tangle}
\object{H\#H^{d}}\step[2.5]\object{H}\\
\O {\lambda'} \step[2]\id\\
\tu {act} \\
\step\object{H}
\end{tangle}\step.
\]
This proved that $\bar{\lambda_{2}}\lambda=\lambda'$ is a injective,
which implies $\bar \lambda \lambda  = id _{H\#
 H^{d}}$.

 (iii) It follows from the simple fact
$\rho (f \# 1) = \lambda (1 \# f)$ for any $f\in U$ (see \cite [
Example 9.4.7] {Mo93}). $\Box$
\begin {Theorem} (Duality Theorem) \label {3.6}  Let $H$ be a  braided  Hopf algebra in $\mathcal{D}=$
$^B_B{\cal YD}$ with invertible antipode, $U$ is a braided Hopf
subalgebra of $H^\circ_{\mathcal{D}}$, the braiding is symmetric on
set $\{B, H, H^{\circ}_{\mathcal{D}}\}$. And let $R$ in $^B_B{\cal
YD}$ be an $U$-comodule algebra such that $R$ is an $H$-module
algebra defined as in section \ref {e1}, $U$ act on $R\#H$ by acting
trivially on $R$ and via $\rightharpoonup $ on $H$. Then

(i) If $U$ satisfy the $RL$-condition with respect to $H$, then
  $$   (R \# H)\# U   \cong R \otimes (H   \# U ) \hbox { \ \ \  as }  \hbox {algebras in  } {\cal D}.    $$

(ii) If $H$ is quantum cocommutative, then
 $$   (R \# H)\# U   \cong R \otimes (H   \# U ) \hbox { \ \ \  as }  \hbox {algebras in } {\cal D}.   $$

\end {Theorem}
{\bf Proof.} It is clear that $U$ is a quasi-dual of $H$ under
evaluation $<,>_{ev} $. $U$ has an invertible antipode since $H$ has
an invertible antipode and $U$ satisfy the $RL$-condition with
respect to $H$ implies that $w$ is a algebra morphism, so
$CRL^{'}$-condition is hold. By Theorem \ref {1.8}, we complete the
proof. $\Box$

This theorem reproduces the main result in \cite {Mo93}.
\begin {Example} \label {3.7} Let $H$ be quantum tenser algebra in ${\cal D}={}^{kG}_{kG}{\cal YD}$
with $\chi_{i}(g_{j})\chi_{j}(g_{i})=1$ for $i,j \in \{1,2, \cdots
\theta\}$ and has invertible antipode, or let $H$ be  a quantum
cocommutative braided Hopf algebra in ${\cal D}={}^B_B{\cal YD}$
with  finite-dimensional commutative and cocomutative $B$ (for
example, $H$ is the unversal enveloping algebra of a Lie
superalgebra). Set $U =H^\circ _{\cal D} = R$. It is clear that $(R,
\phi )$ is a right $U$-comodule algebra with $\phi = \Delta $. By
theorem \ref {3.6} we have
 $$   (R \# U)\# H   \cong R \otimes (U   \# H ) \hbox { \ \ \  as }  \hbox {algebras in } {\cal D}.  $$
\end {Example}


\providecommand{\bysame}{\leavevmode\hbox
to3em{\hrulefill}\thinspace}
\begin{thebibliography}{1}

\bibitem{Bes:cross}
Yu.~N. Bespalov, \emph{Crossed modules and quantum groups in braided
  categories}, Applied Categorical Structures \textbf{5} (1997), no.~2,
  155--204, http://xxx.lanl.gov/abs/q-alg/9510013.

\bibitem{Bes:FRT}
\bysame, \emph{On the braided {FRT}-construction}, J. Nonlinear
Math. Phys.
  \textbf{4} (1997), no.~1-2, 195--205, http://xxx.lanl.gov/abs/q-alg/9510012.

\bibitem{BKLT:int}
Yu.~N. Bespalov, T.~Kerler, V.~V. Lyubashenko, and V.~G. Turaev,
  \emph{Integrals for braided {H}opf algebras}, J. Pure and Appl. Algebra
  \textbf{148} (2000), no.~2, 113--164, Available as
  http://xxx.lanl.gov/abs/q-alg/9709020.

\end{thebibliography}


\begin{thebibliography}{150}
\bibitem {AS02} N. Andruskewisch and H.J.Schneider, Pointed Hopf algebras,
new directions in Hopf algebras, edited by S. Montgomery and H.J.
Schneider, Cambradge University Press, 2002.

\bibitem {AS98}  N. Andruskiewitsch and H. J. Schneider, Lifting of
quantum linear spaces and pointed Hopf algebras of order $p^3$,  J.
Algebra,  {\bf  209}(1998),  658-691.

\bibitem {ATL01} J.Y.Abuhlail, J.Gomez-Torrecillas and F.J.Lobillo, Duality and rational modules in Hopf algebras over
commutative rings, Journal of Algebra, {\bf 240} (2001),
165--184.

\bibitem {BM85}  R. J. Blattner and  S. Montgomery, A
duality theorem for Hopf module algebras,  J. Algebra, {\bf 95}
(1985), 153--172.


\bibitem{BK01}  B. Bakalov,  A. Kirillov,
Lectures on tensor categories and modular functors,
  University Lecture Series Vol. 21, American Mathematical
Society, Providence, 2001.

\bibitem {DZ99} A.Van Daele and Y.Zhang, Galois theory for multiplier Hopf
algebras with integrals, Algebra Representation Theory, {\bf 2}
(1999), 83-106.
\bibitem {DNR01} S.Dascalescu, C.Nastasecu and S. Raianu,
Hopf algebras: an introduction,  Marcel Dekker Inc. , 2001.


\bibitem{Hu05} Yi-Zhi Huang, Vertex operator algebras, the Verlinde conjecture and modular tensor
categories,
 Proc. Nat. Acad. Sci. {\bf 102} (2005) 5352-5356.


\bibitem {He91} M. A. Hennings, Hopf algebras and regular isotopy invariants
for link diagrams, Proc. Cambridge Phil. Soc. ,  {\bf 109 }(1991),
59-77.

\bibitem{HK04} Yi-Zhi Huang, Liang Kong, Open-string vertex algebras, tensor categories and
operads, Commun.Math.Phys. {\bf 250} (2004) 433-471.

\bibitem {Ka97} L. Kauffman , Invariants of links and 3-manifolds via Hopf
algebras, in Geometry and Physics, Marcel Dekker Lecture Notes in
Pure and Appl. Math. {\bf 184} (1997), 471-479.


\bibitem {Ke99}  T.Kerler, Bridged links and tangle presentations of
cobordism categories, Adv. Math. {\bf 141} (1999), 207-- 281.
\bibitem {Mo93}  S. Montgomery, Hopf algebras and their actions on
rings, CBMS Number 82, Published by AMS, 1993.
\bibitem {Ma95b} S. Majid, Foundations of quantum group theory,
Cambridge University Press, 1995.


\bibitem {Ra94}   D.E. Radford,  The trace function and   Hopf
algebras,
 J. Algebra,
{\bf 163} (1994), 583-622.





\bibitem {Ta99} M. Takeuchi, Finite Hopf algebras in braided tensor
categories, J. Pure and Applied Algebras, {\bf 138}(1999), 59-82.

\bibitem {Zh99} Shouchuan Zhang,
Braided Hopf Algebras, Hunan Normal University Press, Changsha,
Second edition. Also in math. RA/0511251.



\bibitem {Zh03} Shouchuan Zhang,
Duality theorem and Drinfeld double in braided tensor categories,
Algebra Colloq. {\bf 10 } (2003)2,  127-134. math.RA/0307255.

\bibitem {ZZH03} Shouchuan Zhang, Yao-Zhong Zhang, Yanying Han, Duality Theorems for
Infinite Braided Hopf Algebras, e-print, math/0309007
\end {thebibliography}
\end{document}